\documentclass[11pt,reqno]{amsart}

\usepackage{amsmath,amsfonts,amssymb,amsthm}
\usepackage[usenames,dvipsnames,svgnames,x11names]{xcolor}
\usepackage{geometry}
\usepackage{enumerate}
\usepackage[pagebackref=true]{hyperref}

\hypersetup{
    colorlinks=true,
    breaklinks=true,
    urlcolor=NavyBlue,
    linkcolor=BrickRed,
    citecolor=ForestGreen,
    bookmarksopen=false,
    filecolor=black
}

\allowdisplaybreaks

\newtheorem{theorem}{Theorem}[section]
\newtheorem{proposition}[theorem]{Proposition}
\newtheorem{lemma}[theorem]{Lemma}
\newtheorem{corollary}[theorem]{Corollary}

\theoremstyle{definition}
\newtheorem{remark}[theorem]{Remark}

\newtheorem{Ithm}{Theorem}

\newtheorem{Iprop}[Ithm]{Proposition}

\newcommand{\HH}{\mathbb H}
\newcommand{\RR}{\mathbb R}
\newcommand{\BB}{\mathbb B}
\newcommand{\Dm}{\mathcal D^{m,2}}
\newcommand{\ISO}{\operatorname{ISO}}
\DeclareMathOperator{\dist}{dist}

\begin{document}

\title[Two-channel compactness for critical GJMS equations]
      {Two-channel global compactness and existence for critical GJMS
       equations on hyperbolic space}
\date{}

\author{Jungang Li}
\address{Jungang Li: Department of Mathematics\\
 University of Science and Technology of China\\
 Hefei, Anhui 230026, China}
\email{jungangli@ustc.edu.cn}

\author{Zhiwei Wang}
\address{Zhiwei Wang: Department of Mathematics\\
 University of Science and Technology of China\\
 Hefei, Anhui 230026, China}
\email{Wzhiwei@mail.ustc.edu.cn}

\begin{abstract}
Let $P_m$ be the GJMS operator of order $2m$ on real hyperbolic space,
$n>2m$.  For the critical equation
\[
    P_m u+a(x)u=|u|^{4m/(n-2m)}u
        \quad\text{on }\HH^n,
\]
where $0\leq a\in L^{n/(2m)}(\HH^n)$, we prove a global compactness
theorem for arbitrary, possibly sign-changing Palais--Smale sequences.
Every such sequence decomposes into a solution of the equation, finitely
many Euclidean profiles concentrating at vanishing scales, finitely many
solutions of the limiting GJMS equation transported to infinity by
divergent hyperbolic isometries, and a remainder converging strongly in the
energy space.  The quadratic form, the $L^{2n/(n-2m)}$ norm and the energy
split along this decomposition.

The main new point is the treatment of the half-space profiles arising at
the conformal boundary of the Poincar\'e ball.  After conformal lifting,
they are precisely finite-energy solutions on $\HH^n$ escaping at a fixed
hyperbolic scale; hence they do not form a third type of profile.  For
sign-changing sequences, Moreau's polar-cone decomposition and positivity
of the limiting Green operators give a sharp two-quantum energy bound and
a compactness range below the two-bubble level.  As an application, we
obtain one solution for sufficiently concentrated nonnegative potentials
and, under an additional smallness condition, a second solution.
\end{abstract}

\maketitle



\section{Introduction}

Let $P_m$ be the GJMS operator of order $2m$ on real hyperbolic space
$\HH^n$, with $n>2m$, and put $q=2n/(n-2m)$.  This paper is concerned
with the critical equation
\[
    P_m u+a(x)u=|u|^{q-2}u,
    \qquad 0\leq a\in L^{n/(2m)}(\HH^n).
\]
Its main object is not a particular family of solutions, but the complete
description of Palais--Smale sequences.  Two independent defects
have to be resolved.  A sequence may concentrate at scales tending to
zero, in which case it sees the Euclidean tangent equation.  It may also
retain a fixed hyperbolic scale and leave every compact set under
diverging isometries.  The first defect is caused by the critical scaling;
the second by the noncompactness of the hyperbolic isometry group.

The principal result of the paper is an intrinsic global compactness
theorem which treats these two defects simultaneously for arbitrary,
possibly sign-changing sequences.  Every Palais--Smale sequence splits
into a solution of the equation with potential, finitely many Euclidean
profiles at vanishing scales, finitely many solutions of the limiting
hyperbolic equation transported to infinity, and a remainder converging
strongly in the energy space.  The quadratic norm, the critical norm and
the energy all decouple.  A point that is specific to the present
argument is that the half-space profiles produced by the ball-model
analysis are identified exactly with the escaping hyperbolic profiles.
Thus conformal infinity creates no additional, coordinate-dependent
channel of loss of compactness.

The equation couples the spectral theory of the higher-order Schr\"odinger
operator $P_m+a$ with the extremal problem for the critical GJMS Sobolev
inequality.  The variable potential is therefore not included merely as a
lower-order perturbation.  The homogeneous critical equation is invariant
under the full isometry group of $\HH^n$; a nonconstant $a$ breaks this
invariance and selects a spatial location.  At the same time,
$L^{n/(2m)}$ is the natural scale-invariant Lebesgue class, since
\[
   \int_{\HH^n}a u^2\,dV
       \leq \|a\|_{L^{n/(2m)}}\|u\|_{L^{2n/(n-2m)}}^2.
\]
The same norm is preserved by the conformal transfer to the Poincar\'e
ball.  In the global decomposition the potential remains in the equation
for the weak limit, but disappears from both the Euclidean blow-up
profiles and the profiles escaping to infinity.  This separates the
inhomogeneous core from two universal limiting equations.  It also makes
localized potentials useful variationally: a concentrated positive
potential can create a nontrivial barycenter--concentration linking on the
topologically trivial space $\HH^n$.  The existence and multiplicity
results below are consequences of this interaction between the potential
and the two compactness defects.

The critical Sobolev equation arose in geometric form in the Yamabe
problem \cite{Yamabe1,Trudinger1,Aubin2,Schoen1,LeeParker1}; the sharp
Euclidean inequality and its extremals were determined by Talenti
\cite{Talenti1}.  The addition of a lower-order spectral term led to the
Br\'ezis--Nirenberg problem \cite{BrezisNirenberg1}.  Topological effects
at the critical exponent were developed by Coron and Bahri--Coron
\cite{Coron,BahriCoron}; see also the account of Br\'ezis \cite{Brezis2}.
Further results on the Br\'ezis--Nirenberg problem may be found in
\cite{SchechterZou1,YueZou}.

For higher-order operators, the critical dimensions and the existence
theory depend more strongly on the order and on the boundary conditions;
see \cite{EdmundsFortunatoJannelli1,PucciSerrin1,Gazzola1,
GazzolaGrunau1}.  Positive solutions of critical polyharmonic Dirichlet
problems were obtained under various hypotheses in
\cite{Grunau2,AlvesdoO}.

The analytic setting for the problem comes from the development of sharp
geometric functional inequalities on hyperbolic space.  Liu
\cite{Liu1} obtained the sharp higher-order Sobolev inequality for the
GJMS energy.  Lu and Yang initiated a systematic higher-order study based
on the Helgason--Fourier transform and the hyperbolic Plancherel formula:
their results include a sharp Hardy--Adams inequality for the
bi-Laplacian \cite{LuYangHardyAdams}, higher-order
Hardy--Sobolev--Maz'ya inequalities \cite{LuYang1}, and precise Green
functions and sharp constants for Paneitz and GJMS operators
\cite{LuYang2}.  The Helgason--Fourier calculus is particularly well
suited to the factorization of $P_m$ as a polynomial in the hyperbolic
Laplacian; it reveals the spectral Hardy terms and gives positive Green
kernels even though a differential maximum principle is generally
unavailable in higher order.

More recently, Lu and Tao \cite{LuTao1} proved the existence and symmetry
of extremals for the critical higher-order Hardy--Sobolev--Maz'ya
inequality when $n\geq 2m+2$.  In the equivalent hyperbolic formulation,
their result gives extremals for the corresponding Poincar\'e--Sobolev
inequality and positive radial solutions of the higher-order
Br\'ezis--Nirenberg equation at the bottom of the $L^2$-spectrum of $P_m$.
Their proof combines the Helgason--Fourier transform and positive Green
kernels with hyperbolic rearrangement and concentration--compactness.

Building on these inequalities, Li, Lu and Yang \cite{LiLuYang1}
established the higher-order Br\'ezis--Nirenberg theory on bounded domains
and on the whole of $\HH^n$, including existence, nonexistence and
symmetry for the constant spectral perturbation
\[
    P_m u-\lambda u=|u|^{q-2}u.
\]
The present work concerns a variable coefficient at the critical form
scale.  Its purpose is to describe every Palais--Smale sequence rather
than to establish the existence or symmetry of a distinguished positive
solution.

Classification results also enter at the Euclidean end of the
decomposition.  Positive finite-energy solutions of the critical
polyharmonic equation on $\RR^n$ are the standard Sobolev bubbles.  This
was proved by Lin \cite{Lin} in the fourth-order case and by Wei and Xu
\cite{WeiXu} for every integer order; the integral-equation
classification of Chen, Li and Ou \cite{ChenLiOu} gives a complementary
route.  In particular, the super-polyharmonic property established by
Wei and Xu is the basic PDE mechanism that permits passage to the positive
integral formulation.  These results identify the positive Euclidean
energy quantum.  Our decomposition itself does not assume positivity and
therefore retains general sign-changing Euclidean and hyperbolic
profiles.

Classical symmetry arguments for second-order semilinear equations rely
on the maximum principle and the method of moving planes
\cite{GidasNiNirenberg1,BerestyckiNirenberg1}.  Unique continuation for
certain second- and fourth-order equations was established by Pederson
\cite{Pederson}.  These tools are not available uniformly at arbitrary
order or for sign-changing profiles, and no classification of such
profiles is assumed here.

The second-order equation on $\HH^n$ has also been studied extensively.
Positive and radial solutions were considered in
\cite{ManciniSandeep1,BandleKabeya1}; existence, uniqueness and solution
gaps for the hyperbolic Br\'ezis--Nirenberg problem were studied in
\cite{Stapelkamp1,Benguria1}.  Nondegeneracy and sign-changing solutions
were treated in \cite{GangulySandeep1,GangulySandeep2}.

For $m=1$, Bhakta and Sandeep \cite{BhaktaSandeep1} already exhibited
the two compactness channels for the critical Poincar\'e--Sobolev
equation.  Bhakta, Ganguly, Gupta and Sahoo
\cite{BhaktaGangulyGuptaSahoo} subsequently obtained a two-profile global
compactness result for a second-order equation with a nonautonomous
critical coefficient and a forcing term.  Profile decompositions in
Sobolev spaces on noncompact manifolds were developed by Sandeep and
Tintarev \cite{SandeepTintarev}.  The higher-order situation considered
here is different in two respects: the perturbation is the form-critical
linear potential $a(x)u$, and arbitrary Palais--Smale sequences may change
sign.  After conformal transfer we use the bounded-domain polyharmonic
decomposition of Mazumdar \cite{Mazumdar} to extract whole-space and
half-space profiles.  The latter are then identified, together with their
parameters and energies, as fixed-scale solutions escaping under
hyperbolic isometries.

\subsection{The equation and its conformal models}

Fix $n>2m$, and put
\[
    q:=\frac{2n}{n-2m},
    \qquad
    a_0:=\frac{n-2m}{2}.
\]
Let $P_m$ be the GJMS operator of order $2m$ on $\HH^n$
\cite{GrahamJenneMasonSparling1,Juhl1,Juhl2}.  We study
\begin{equation}\label{Hyperbolic equation}
    \begin{cases}
        P_m u+a(x)u=|u|^{q-2}u \quad&\text{on }\HH^n,\\
        u\in\Dm(\HH^n),
    \end{cases}
\end{equation}
where
\[
    0\leq a\in L^{n/(2m)}(\HH^n).
\]
Here $\Dm(\HH^n)$ is the completion of $C_c^\infty(\HH^n)$ in the norm
induced by the quadratic form of $P_m$.  For an open set
$D\subset\RR^n$, we write $\Dm(D)$ for the completion of
$C_c^\infty(D)$ in $\|\nabla^m\cdot\|_{L^2}$.  We also fix a base point
$o\in\HH^n$.

Let $\iota$ identify the Poincar\'e ball with $\HH^n$, and set
\[
    \Phi(x):=\frac{2}{1-|x|^2},
    \qquad
    \mathcal Cu:=\Phi^{a_0}(u\circ\iota).
\]
Conformal covariance makes $\mathcal C$ an isometry of the critical
quadratic forms and transforms \eqref{Hyperbolic equation} into
\begin{equation}\label{Euclidean equation}
    \begin{cases}
        (-\Delta)^m v+(a\circ\iota)(x)\Phi(x)^{2m}v=|v|^{q-2}v
            \quad&\text{on }\BB^n,\\
        v\in W_0^{m,2}(\BB^n).
    \end{cases}
\end{equation}
The exponent $n/(2m)$ is exactly the one for which this transformation
preserves the norm of the potential:
\[
    \int_{\BB^n}
       \bigl[(a\circ\iota)(x)\Phi(x)^{2m}\bigr]^{n/(2m)}\,dx
       =
    \int_{\HH^n}a^{n/(2m)}\,dV.
\]
Thus the potential is at the borderline dictated by the critical
Sobolev scaling.

We shall use three limiting equations.  At hyperbolic infinity one has
\begin{equation}\label{hyperbolic infinity problem}
    \begin{cases}
        P_m u=|u|^{q-2}u \quad&\text{on }\HH^n,\\
        u\in\Dm(\HH^n),
    \end{cases}
\end{equation}
while the corresponding Euclidean problems are
\begin{equation}\label{Euclidean unit ball problem}
    \begin{cases}
        (-\Delta)^m v=|v|^{q-2}v \quad&\text{in }\BB^n,\\
        v\in W_0^{m,2}(\BB^n),
    \end{cases}
\end{equation}
and
\begin{equation}\label{Euclidean Whole Space problem}
    \begin{cases}
        (-\Delta)^m v=|v|^{q-2}v \quad&\text{in }\RR^n,\\
        v\in\Dm(\RR^n).
    \end{cases}
\end{equation}
The energy functionals associated with
\eqref{Hyperbolic equation}, \eqref{hyperbolic infinity problem} and
\eqref{Euclidean Whole Space problem} are
\[
    \begin{aligned}
        I(u)&:=\frac12\int_{\HH^n}uP_m u\,dV
              +\frac12\int_{\HH^n}a(x)u^2\,dV
              -\frac1q\int_{\HH^n}|u|^q\,dV,\\
        I_\infty(u)&:=\frac12\int_{\HH^n}uP_m u\,dV
              -\frac1q\int_{\HH^n}|u|^q\,dV,\\
        J(v)&:=\frac12\int_{\RR^n}v(-\Delta)^m v\,dx
              -\frac1q\int_{\RR^n}|v|^q\,dx.
    \end{aligned}
\]
We denote by $S_{m,n}$ the sharp Sobolev constant of $P_m$ on $\HH^n$.
By a theorem of Liu \cite{Liu1}, it agrees with the Euclidean sharp
constant \cite{Lieb,Swanson}.  The remaining geometric and
functional-analytic facts are collected in
Section~\ref{sec:preliminaries}.

\subsection{The two channels of loss of compactness}

Let $\{u_k\}$ be a Palais--Smale sequence for $I$, transfer it to the ball
by $\mathcal C$, and subtract its weak limit.  The term involving the
critical potential is compact in the dual Sobolev topology.  The residual
sequence is therefore Palais--Smale for the pure polyharmonic functional
on $\BB^n$, to which the bounded-domain decomposition applies.  Its
interior profiles solve \eqref{Euclidean Whole Space problem}.  Its
boundary profiles, however, solve the critical Dirichlet problem on a
half-space.

On a general compact manifold with boundary, a half-space profile is a
genuine boundary phenomenon.  In the present problem the boundary is
introduced only by conformal compactification.  The following result
identifies what the half-space alternative means intrinsically on
$\HH^n$.

\begin{Ithm}\label{intro boundary identification}
Endow $\RR^n_+=\{(z,t):t>0\}$ with the hyperbolic metric
$g_+=t^{-2}(dz^2+dt^2)$.  If $W\in\Dm(\RR^n_+)$ solves
\[
    (-\Delta)^mW=|W|^{q-2}W \ \text{ in }\RR^n_+,
    \qquad
    \partial_\nu^\ell W=0 \ \text{ on }\partial\RR^n_+,
    \quad 0\leq\ell\leq m-1,
\]
then its conformal lift $V(z,t):=t^{a_0}W(z,t)$ lies in $\Dm(\HH^n)$,
solves \eqref{hyperbolic infinity problem}, and has the same quadratic
norm, $L^q$ norm and critical energy as $W$.  Under this correspondence
the Euclidean boundary translations and dilations
$(z,t)\mapsto(\zeta+\varepsilon z,\varepsilon t)$ become hyperbolic
isometries, and every boundary bubble of the ball-model decomposition is
of the form $V\circ\eta_k^{-1}$ with $\eta_k\in\ISO(\HH^n)$ and
$\eta_k(o)\to\infty$, up to a strongly vanishing error and a fixed
tangential orthogonal reparametrization of $W$.
\end{Ithm}

The complete statement, including the normalization of the boundary
charts, is Theorem~\ref{half space isometry identification}.  Its proof is
geometric and independent of the order $m$.  It uses the upper
half-space model to track the conformal multiplier, the tangential
translations and the dilations exactly.  In particular, it neither
classifies the half-space solutions nor attempts to rule them out by a
sign or Liouville argument.

In the second-order argument of Bhakta and Sandeep
\cite{BhaktaSandeep1}, the escaping profiles are first detected by
recentring with hyperbolic isometries, after which the remaining critical
concentration is localized and treated by the classical bounded-domain
theory.  For higher-order equations the local polyharmonic decomposition
also permits half-space profiles.  It is used here only for the local
extraction.  What must then be proved is that the half-space parameters,
profiles and energies have an intrinsic hyperbolic meaning.
Theorem~\ref{intro boundary identification} provides exactly this
conversion, without a positivity assumption and without classifying or
excluding the half-space solution.

\begin{Ithm}\label{Global compactness}
Let $a\geq0$ belong to $L^{n/(2m)}(\HH^n)$, and let
$\{u_k\}\subset\Dm(\HH^n)$ be a Palais--Smale sequence for $I$ at level
$d$.  After passing to a subsequence, there exist a solution $u$ of
\eqref{Hyperbolic equation} and integers $N_{\rm E},N_{\rm H}\geq0$
together with the following data.
\begin{enumerate}[(a)]
    \item For $1\leq\alpha\leq N_{\rm E}$, a nonzero solution
    $U^\alpha\in\Dm(\RR^n)$ of \eqref{Euclidean Whole Space problem},
    concentration parameters $\varepsilon_k^\alpha\to0$ and centres
    $x_k^\alpha\in\BB^n$ with
    \[
        \frac{\dist(x_k^\alpha,\partial\BB^n)}
             {\varepsilon_k^\alpha}\longrightarrow\infty,
    \]
    and cut-off functions $\chi_k^\alpha\in C_c^\infty(\BB^n)$ equal to
    $1$ near $x_k^\alpha$.  The associated concentrating bubble in the
    ball model is
    \begin{equation}\label{intro euclidean bubble}
        \mathcal B_k^\alpha(x)
           :=\chi_k^\alpha(x)\,
             (\varepsilon_k^\alpha)^{-a_0}\,
             U^\alpha\!\left(
               \frac{x-x_k^\alpha}{\varepsilon_k^\alpha}\right).
    \end{equation}

    \item For $1\leq\beta\leq N_{\rm H}$, a nonzero solution
    $V^\beta\in\Dm(\HH^n)$ of \eqref{hyperbolic infinity problem} and
    isometries $\eta_k^\beta\in\ISO(\HH^n)$ with
    $\eta_k^\beta(o)\to\infty$.
\end{enumerate}
With this notation,
\begin{equation}\label{intro global decomposition}
    u_k=u+\sum_{\alpha=1}^{N_{\rm E}}
                \mathcal C^{-1}\mathcal B_k^\alpha
          +\sum_{\beta=1}^{N_{\rm H}}
                V^\beta\circ(\eta_k^\beta)^{-1}+r_k,
    \qquad
    r_k\to0\quad\text{in }\Dm(\HH^n).
\end{equation}
The parameters of distinct Euclidean bubbles are asymptotically
orthogonal, the hyperbolic isometries are mutually divergent, and the
quadratic and $L^q$ norms decouple accordingly.  In particular,
\begin{equation}\label{intro global energy splitting}
    d=I(u)+\sum_{\alpha=1}^{N_{\rm E}}J(U^\alpha)
          +\sum_{\beta=1}^{N_{\rm H}}I_\infty(V^\beta).
\end{equation}
\end{Ithm}

Theorem~\ref{Global compactness} is the principal result of the paper.
The two lists of profiles describe different geometric limits.  The
$U^\alpha$ arise after magnifying a vanishing scale and therefore see the
Euclidean tangent equation.  The $V^\beta$ retain their hyperbolic scale
and disappear only because the isometries leave every compact subset.
Theorem~\ref{intro boundary identification} shows, in particular, that a
third boundary channel is not needed: the apparent boundary profiles of
the ball model are exactly the second, intrinsic channel.

\subsection{The compactness threshold}

The decomposition \eqref{intro global decomposition} reduces compactness
to lower bounds for the energies of the nonzero profiles.  In the
second-order problem, the positive and negative parts of a sign-changing
solution can be tested separately, and each carries one Sobolev quantum.
This argument is unavailable in $W^{m,2}$ when $m\geq2$, since
$u^+$ and $u^-$ need not belong to the energy space.  We replace it by an
order-theoretic decomposition.

Let $(X,\mu)$ be a measure space, let $2<q<\infty$, and let $H$ be a real
Hilbert space with inner product $\mathfrak a$, continuously embedded in
$L^q(X,\mu)$.  Write
\[
    \mathcal K:=\{u\in H:u\geq0\ \mu\text{-a.e.}\},
    \qquad
    \mathcal K^\circ:=\{w\in H:\mathfrak a(w,v)\leq0
        \ \text{ for all } v\in\mathcal K\},
\]
and set
\[
    S_H:=\inf_{u\in H\setminus\{0\}}
        \frac{\mathfrak a(u,u)}{\|u\|_{L^q}^2}>0.
\]

\begin{Ithm}\label{intro nodal principle}
Suppose $\mathcal K^\circ\subset-\mathcal K$, which happens if and only if
the inverse of the Riesz map of $\mathfrak a$ is positivity preserving.
Then every sign-changing critical point $u$ of
\[
    \mathcal E(u)=\tfrac12\mathfrak a(u,u)-\tfrac1q\|u\|_{L^q}^q
\]
satisfies
\[
    \mathcal E(u)\geq
    2\Bigl(\frac12-\frac1q\Bigr)S_H^{q/(q-2)}.
\]
\end{Ithm}

The decomposition of a higher-order Sobolev function with respect to
polar cones is classical; see
\cite[Section~3.1.2]{GazzolaGrunauSweers1}.  Moreau's theorem
\cite{Moreau} splits $u$ orthogonally along $\mathcal K$ and
$\mathcal K^\circ$.  Under the inclusion in
Theorem~\ref{intro nodal principle}, the two components have opposite
pointwise signs, and each contributes at least one Sobolev quantum.  The
factor $2$ is optimal in this abstract setting.  The theorem applies to
$(-\Delta)^m$ on $\RR^n$, whose Riesz kernel is positive, and to $P_m$ on
$\HH^n$, using the positive Green kernel found by Lu and Yang
\cite{LuYang2}.

For a hyperbolic profile there is one further piece of information.  A
nonexistence result on the ball
\cite[Theorem~7.34]{GazzolaGrunauSweers1} rules out a solution of one
sign to \eqref{hyperbolic infinity problem}.  Such nonexistence results
originate in the identity of Pohozaev \cite{Pohozaev1} and its general
variational extension due to Pucci and Serrin \cite{PucciSerrin2}.  Hence
every nonzero hyperbolic profile is sign-changing and satisfies
\[
    I_\infty(V)\geq\frac{2m}{n}S_{m,n}^{n/(2m)}.
\]
Combined with Theorem~\ref{Global compactness}, this gives the
Palais--Smale condition for the normalized constrained functional in the
window
\[
    \bigl(S_{m,n},2^{2m/n}S_{m,n}\bigr).
\]

\subsection{Existence for concentrated critical potentials}

This application illustrates the role of the variable potential.  The
homogeneous equation cannot distinguish one point of $\HH^n$ from
another, and a variational family may disappear either by concentrating
or by moving to infinity.  A concentrated positive potential breaks this
transitivity and creates a nontrivial linking in the barycenter and scale
variables.  The global compactness theorem prevents the corresponding
minimax sequences from losing compactness through either of the two
channels.

Variational methods for nonlinear scalar-field equations were developed
systematically by Berestycki and Lions
\cite{BerestyckiLions1,BerestyckiLions2}.  The construction used here is a
higher-order hyperbolic version of the potential mechanism of Benci and
Cerami \cite{BenciCerami}.  In the second-order Euclidean problem, a
potential in the critical class can replace the topology of the domain in
the existence argument.  Passaseo \cite{Passaseo} subsequently considered
a potential concentrated at one point and obtained one or two solutions,
according to the size of its critical norm.

For the present problem, Theorem~\ref{abstract two level criterion}
provides two critical points of a $C^1$ functional on a Finsler manifold
once a barycenter--concentration map has nonzero degree and the associated
minimax levels lie in the compactness window above.  Its application gives
the following result.

\begin{Ithm}\label{Main Theorem}
Let $n>2m$ and fix $x_0\in\BB^n$.  Let
$\bar\alpha\in L^{n/(2m)}(\BB^n)$ and
$\alpha\in L^{n/(2m)}(\RR^n)$ be nonnegative, with
$\alpha\not\equiv0$.  Consider \eqref{Hyperbolic equation} with potential
determined by
\[
    \Phi(x)^{2m}(a\circ\iota)(x)
       =\bar\alpha(x)+\lambda^{2m}
          \alpha\bigl(\lambda(x-x_0)\bigr),
\]
and, on the $L^q$-unit sphere of $W_0^{m,2}(\BB^n)$, set
\[
    f_\lambda(v):=
       \int_{\BB^n}\left(
          |\nabla^m v|^2+
          \bigl[\bar\alpha(x)+\lambda^{2m}
             \alpha\bigl(\lambda(x-x_0)\bigr)\bigr]v^2
       \right)dx.
\]
\begin{enumerate}[(i)]
    \item For every $\epsilon>0$ there is $\Lambda(\epsilon)>0$ such that
    for every $\lambda\geq\Lambda(\epsilon)$ the problem
    \eqref{Hyperbolic equation} has a nontrivial weak solution $u_\lambda$
    whose $L^q$-normalized representative $v_\lambda$ in the ball model
    satisfies
    \[
        S_{m,n}<f_\lambda(v_\lambda)<S_{m,n}+\epsilon.
    \]
    One may therefore select a family of such solutions, indexed by all
    sufficiently large $\lambda$, along which
    $f_\lambda(v_\lambda)\to S_{m,n}$.

    \item If moreover
    \[
        \|\alpha\|_{L^{n/(2m)}(\RR^n)}
           <S_{m,n}\bigl(2^{2m/n}-1\bigr),
    \]
    then $\Lambda(\epsilon)$ may be chosen so that for every
    $\lambda\geq\Lambda(\epsilon)$ there is a second nontrivial weak
    solution $\widehat u_\lambda$ whose normalized ball representative
    $\widehat v_\lambda$ satisfies
    \[
        f_\lambda(v_\lambda)<f_\lambda(\widehat v_\lambda)
           <2^{2m/n}S_{m,n}.
    \]
    In particular the two solutions are distinct even modulo sign.
\end{enumerate}
\end{Ithm}

\subsection{Why no positivity conclusion is made}
\label{sec:signobstruction}

Theorem~\ref{Main Theorem} asserts nothing about the sign of the solutions
it produces.  This is not merely a limitation of the minimax
construction.  For higher-order operators, the hypotheses of the theorem
do not ensure inverse positivity of the perturbed operator throughout the
admissible class.

\begin{Iprop}\label{higher-order sign obstruction}
Let $m>1$.  There exist $\bar\alpha$ and $\alpha$ satisfying all the
hypotheses of Theorem~\ref{Main Theorem}, including the smallness
condition in part \textup{(ii)}, for which the operator
\[
    (-\Delta)^m+\bar\alpha(x)+\lambda^{2m}
        \alpha\bigl(\lambda(x-x_0)\bigr)
\]
on $W^{m,2}_0(\BB^n)$ has an inverse that fails to be positivity
preserving, for every $\lambda>0$.
\end{Iprop}

\begin{proof}
By Corollary~6.4 of \cite{GrunauSweersPositivity}, there is a constant
$A_c\geq0$ such that $((-\Delta)^m+V)^{-1}$ is not positivity preserving
on the ball whenever a continuous potential exceeds $A_c$ throughout the
ball.  Take $\bar\alpha\equiv A$ with $A>A_c$, and let $\alpha$ be any
nonzero nonnegative function in $C_c^\infty(\RR^n)$.  The amplitude of
$\alpha$ may be taken arbitrarily small, in particular small enough for
the hypothesis of part~\textup{(ii)}, whereas
\[
    \bar\alpha(x)+\lambda^{2m}
        \alpha\bigl(\lambda(x-x_0)\bigr)>A_c
\]
for every $\lambda>0$ and every $x\in\BB^n$.
\end{proof}

Proposition~\ref{higher-order sign obstruction} does not say that positive
solutions fail to exist.  It says that their positivity cannot be deduced
uniformly from an inverse-positivity argument under the hypotheses of
Theorem~\ref{Main Theorem}.  Inverse positivity is used in this paper only
for the unperturbed limiting operators, through
Theorem~\ref{intro nodal principle}, where the relevant Green kernels are
positive.

\subsection{Organization of the paper}

Section~\ref{sec:preliminaries} recalls the conformal models, the GJMS
energy space and the sharp Sobolev inequality.
Section~\ref{sec:struwe} states the bounded-domain polyharmonic
decomposition and proves Theorem~\ref{intro boundary identification}.
Section~\ref{sec:globalcompactness} proves
Theorem~\ref{Global compactness}.  Section~\ref{sec:moreau} develops the
abstract Moreau principle, proves
Theorem~\ref{intro nodal principle} and derives the Palais--Smale window.
Finally, Section~\ref{sec:existence} establishes the two-level
barycenter--concentration criterion and applies it to prove
Theorem~\ref{Main Theorem}.

\section{Notation and preliminaries}\label{sec:preliminaries}

Throughout the paper, $m\in\mathbb N$, $n>2m$, and
\[
    q:=\frac{2n}{n-2m},
    \qquad a_0:=\frac{n-2m}{2}.
\]
We collect here the geometric and variational facts used later on,
referring to \cite{GrahamJenneMasonSparling1,Juhl1,Juhl2} for the GJMS
operators and to \cite{Liu1,LuYang1,LuYang2} for the sharp inequalities.

\subsection{The half-space and ball models}

The upper-half-space model is
\[
    \mathbb R^n_+
       =\{(z,t)\in\mathbb R^{n-1}\times\mathbb R:t>0\},
    \qquad
    g_+=t^{-2}(dz^2+dt^2),
    \qquad
    dV_{g_+}=t^{-n}\,dz\,dt.
\]
The ball model is
\[
    \mathbb B^n=\{x\in\mathbb R^n:|x|<1\},
    \qquad
    g_{\mathbb B}=\Phi(x)^2g_{\rm Eucl},
    \qquad
    \Phi(x):=\frac{2}{1-|x|^2},
\]
with $dV_{g_{\mathbb B}}=\Phi^n\,dx$.  A Cayley transformation is an
isometry between these two models.  We write
\[
    \iota:(\mathbb B^n,g_{\mathbb B})\longrightarrow\mathbb H^n
\]
for the ball-model identification.

\subsection{The GJMS energy space}\label{GJMS eigenvalues}

On $\mathbb H^n$, the order-$2m$ GJMS operator factors as
\begin{equation}\label{GJMS factorization}
    P_m=\prod_{\ell=1}^m\bigl(P_1+\ell(\ell-1)\bigr),
    \qquad
    P_1=-\Delta_{\mathbb H^n}-\frac{n(n-2)}4.
\end{equation}
Its quadratic form is positive on $C_c^\infty(\mathbb H^n)$.  We set
\[
    \|u\|_{\mathcal D^{m,2}(\mathbb H^n)}^2
       :=\int_{\mathbb H^n}uP_m u\,dV
\]
and define $\mathcal D^{m,2}(\mathbb H^n)$ as the completion of
$C_c^\infty(\mathbb H^n)$ in this norm.

For Euclidean domains, we use the convention
\[
   \int_\Omega|\nabla^m v|^2\,dx
   :=
   \begin{cases}
      \displaystyle\int_\Omega|\Delta^{m/2}v|^2\,dx,
          &m\ \text{even},\\[4pt]
      \displaystyle\int_\Omega
           |\nabla\Delta^{(m-1)/2}v|^2\,dx,
          &m\ \text{odd}.
   \end{cases}
\]
This is the quadratic form of $(-\Delta)^m$ on
$W_0^{m,2}(\Omega)$.

The conformal covariance of $P_m$ takes the form
\begin{equation}\label{preliminary conformal covariance}
    (P_m u)\circ\iota
      =\Phi^{-\frac{n+2m}{2}}(-\Delta)^m
          \bigl[\Phi^{a_0}(u\circ\iota)\bigr].
\end{equation}
Consequently, the transfer
\begin{equation}\label{preliminary conformal transfer}
    \mathcal Cu:=\Phi^{a_0}(u\circ\iota)
\end{equation}
extends by density to an isometric isomorphism
\[
    \mathcal C:\mathcal D^{m,2}(\mathbb H^n)
       \longrightarrow W_0^{m,2}(\mathbb B^n),
\]
and
\begin{align}
    \int_{\mathbb H^n}uP_m u\,dV
       &=\int_{\mathbb B^n}|\nabla^m\mathcal Cu|^2\,dx,
       \label{preliminary quadratic identity}\\
    \int_{\mathbb H^n}|u|^q\,dV
       &=\int_{\mathbb B^n}|\mathcal Cu|^q\,dx.
       \label{preliminary Lq identity}
\end{align}

\subsection{The sharp Sobolev constant}

The sharp GJMS Sobolev inequality \cite{Liu1} is
\begin{equation}\label{sharp GJMS Sobolev inequality}
    S_{m,n}
       \left(\int_{\mathbb H^n}|u|^q\,dV\right)^{2/q}
       \leq\int_{\mathbb H^n}uP_m u\,dV,
       \qquad u\in\mathcal D^{m,2}(\mathbb H^n),
\end{equation}
where $S_{m,n}$ is the Euclidean sharp constant.  This constant is not
attained by a nonzero function in
$\mathcal D^{m,2}(\mathbb H^n)$.  Indeed, an extremal would be carried by
\eqref{preliminary conformal transfer} to an extremal in
$W_0^{m,2}(\mathbb B^n)$; extension by zero would then be a whole-space
Sobolev extremal vanishing on an open set, contrary to the
equality-case classification for the higher-order sharp Sobolev inequality
\cite{Lieb,Swanson}.

We write $S=S_{m,n}$ whenever the indices play no role.  The one-bubble
energy and the corresponding normalized quadratic level are, respectively,
\[
    \frac mnS_{m,n}^{n/(2m)}
    \qquad\text{and}\qquad
    S_{m,n}.
\]

\subsection{Palais--Smale sequences}

Let $\mathcal X$ be a $C^1$ Banach or Finsler manifold and
$F\in C^1(\mathcal X)$.  A sequence $\{u_k\}\subset\mathcal X$ is a
Palais--Smale sequence at level $c$ if
\[
    F(u_k)\longrightarrow c,
    \qquad
    \|dF(u_k)\|_{T_{u_k}^*\mathcal X}\longrightarrow0.
\]
We say that $F$ satisfies the Palais--Smale condition at $c$ if every
such sequence has a strongly convergent subsequence.

\section{Polyharmonic profile decomposition on the ball}
\label{sec:struwe}

The argument of Section~\ref{sec:globalcompactness} is carried out after
the problem has been transferred conformally to the Euclidean unit ball,
where concentration may occur both in the interior and at the boundary.
We therefore work with the bounded-domain polyharmonic decomposition of
Mazumdar \cite{Mazumdar}, the polyharmonic analogue of Struwe's theorem
\cite{Struwe}, recalled below as
Theorem~\ref{Euclidean Struwe Decomposition}.  What is new here is
Theorem~\ref{half space isometry identification}.

For an open set $D\subset\mathbb R^n$, write
\[
    \mathcal D^{m,2}(D)
    :=\overline{C_c^\infty(D)}^{\,\|\nabla^m\cdot\|_2}
\]
and set
\[
    E_D(v):=\frac12\int_D|\nabla^m v|^2\,dx
             -\frac1q\int_D|v|^q\,dx,
    \qquad q=\frac{2n}{n-2m}.
\]
We use $\mathbb R^n_+:=\{(z,t)\in\mathbb R^{n-1}\times\mathbb R:t>0\}$.
Mazumdar formulates the boundary alternative on $\mathbb R^n_-$; throughout
we compose his boundary charts with the fixed reflection in the last
coordinate.  The boundary energy space in that theorem is the clamped
half-space realization used below, with the trace conditions in
\eqref{half space profile equation}.

\begin{theorem}[Mazumdar's decomposition on a bounded domain]
\label{Euclidean Struwe Decomposition}
Let $m\in\mathbb N$, $n>2m$, let
$\Omega\subset\mathbb R^n$ be a smooth bounded domain, and let
$\{w_k\}\subset W^{m,2}_0(\Omega)$ be a Palais--Smale sequence for
$E_\Omega$.  Then, after passing to a subsequence, there are a weak
solution $w^0\in W^{m,2}_0(\Omega)$, finitely many nonzero profiles
$W^1,\ldots,W^N$, scales $\varepsilon_k^j\to0$, and centres
$c_k^j\in\overline\Omega$ such that
\begin{equation}\label{bounded profile decomposition}
    w_k=w^0+\sum_{j=1}^N\mathcal B_k^j+r_k,
    \qquad r_k\longrightarrow0
    \quad\text{strongly in }W^{m,2}_0(\Omega).
\end{equation}
Each profile is of precisely one of the following two types.
\begin{enumerate}[(i)]
    \item \emph{Interior profile.}  One has $c_k^j\in\Omega$ and
    \[
        \frac{\dist(c_k^j,\partial\Omega)}
             {\varepsilon_k^j}\longrightarrow\infty,
    \]
    while $W^j\in\mathcal D^{m,2}(\mathbb R^n)$ solves
    \[
        (-\Delta)^mW^j=|W^j|^{q-2}W^j
        \quad\text{in }\mathbb R^n.
    \]
    Up to a term converging strongly to zero, $\mathcal B_k^j$ is a
    cut-off of
    \[
        (\varepsilon_k^j)^{-\frac{n-2m}{2}}
        W^j\!\left(\frac{x-c_k^j}{\varepsilon_k^j}\right).
    \]

    \item \emph{Boundary profile.}  The centres may be chosen on
    $\partial\Omega$, and $W^j\in\mathcal D^{m,2}(\mathbb R^n_+)$ solves
    \begin{equation}\label{half space profile equation}
        \begin{cases}
        (-\Delta)^mW^j=|W^j|^{q-2}W^j&\text{in }\mathbb R^n_+,\\
        \partial_\nu^\ell W^j=0&\text{on }\partial\mathbb R^n_+,
            \quad 0\leq\ell\leq m-1.
        \end{cases}
    \end{equation}
    The boundary conditions are understood in the trace, equivalently
    weak, sense.  If $\mathcal T_c$ is a smooth boundary-flattening chart whose
    differential at the origin is an isometry, then, again up to a
    strongly vanishing term, $\mathcal B_k^j$ is a cut-off of
    \[
        (\varepsilon_k^j)^{-\frac{n-2m}{2}}
        W^j\!\left(
        \frac{\mathcal T_{c_k^j}^{-1}(x)}{\varepsilon_k^j}
        \right).
    \]
\end{enumerate}
The parameters are asymptotically orthogonal: for $i\neq j$,
\begin{equation}\label{Mazumdar parameter orthogonality}
    \frac{\varepsilon_k^i}{\varepsilon_k^j}
    +\frac{\varepsilon_k^j}{\varepsilon_k^i}
    +\frac{|c_k^i-c_k^j|^2}
          {\varepsilon_k^i\varepsilon_k^j}
    \longrightarrow\infty.
\end{equation}
Moreover, with $D_j=\mathbb R^n$ for an interior profile and
$D_j=\mathbb R^n_+$ for a boundary profile,
\begin{align}
    E_\Omega(w_k)
      &=E_\Omega(w^0)+\sum_{j=1}^NE_{D_j}(W^j)+o(1),
      \label{bounded energy splitting}\\
    \int_\Omega|\nabla^m w_k|^2\,dx
      &=\int_\Omega|\nabla^m w^0|^2\,dx
        +\sum_{j=1}^N\int_{D_j}|\nabla^mW^j|^2\,dx+o(1),
      \label{bounded norm splitting}\\
    \int_\Omega|w_k|^q\,dx
      &=\int_\Omega|w^0|^q\,dx
        +\sum_{j=1}^N\int_{D_j}|W^j|^q\,dx+o(1).
      \label{bounded Lq splitting}
\end{align}
\end{theorem}

The decomposition itself is quoted, not proved here.  We supply only the
splitting of the quadratic and $L^q$ norms, which we shall use in the form
stated.

\begin{proof}[Attribution and proof of the splitting identities]
The decomposition is the specialization of \cite[Theorem~1.1]{Mazumdar} to
the polyharmonic operator on a smooth bounded Euclidean domain.  The two
alternatives are those in \cite[Definition~2.1]{Mazumdar}; independence
of the cut-off and boundary chart, modulo a strongly convergent error,
is \cite[Proposition~2.1]{Mazumdar}.  The parameter separation is proved
in \cite[Proposition~5.1]{Mazumdar}.  In particular, the quadratic
cross terms between distinct profiles tend to zero.  Write
$B_k^i,B_k^j$ for two cut-off bubbles.  First approximate their
profiles in the homogeneous space by compactly supported smooth functions.
After rescaling one of the two bubbles, the parameter separation gives
\[
    \int_\Omega\nabla^m B_k^i\cdot\nabla^m B_k^j\,dx=o(1),
    \qquad
    \int_\Omega\bigl(|B_k^i+B_k^j|^q
         -|B_k^i|^q-|B_k^j|^q\bigr)\,dx=o(1).
\]
The uniform boundedness of the bubble operators then removes the smooth
approximation.  The same argument separates the weak limit from every
bubble.  Iterated use of the Br\'ezis--Lieb lemma \cite{BrezisLieb1},
together with the strong remainder, gives the displayed $L^q$ splitting;
the first estimate gives the quadratic splitting.  The Palais--Smale
identities then yield the energy splitting.
\end{proof}

We turn to the geometric description of the boundary profiles.  This is the
step that converts the two Euclidean alternatives of Mazumdar's theorem into
the two hyperbolic alternatives of Theorem~\ref{Global compactness}.

A preliminary normalization of the charts is required.  A Cayley chart
carries a Euclidean conformal factor which cannot be discarded outright; it
must first be absorbed into the scale, and only then may one pass to the
profile limit.

\begin{lemma}[Normalization of boundary charts]
\label{Cayley chart normalization}
Let $W\in\mathcal D^{m,2}(\mathbb R^n_+)$, and let
$\mathcal B_k^{\mathcal T}[W]$ be a boundary bubble, with centres converging
to a fixed point of $\partial\mathbb B^n$, scales $\varepsilon_k\to0$, and
an admissible flattening chart in the sense of
Theorem~\ref{Euclidean Struwe Decomposition}.  Put $a_0=(n-2m)/2$.
Then, after passing to a subsequence and changing the representative of the
bubble by a strongly vanishing term, there are
\[
    \frac{\varepsilon_k'}{\varepsilon_k}\longrightarrow1,
    \qquad
    O_k\in O(n-1),\quad O_k\longrightarrow O_\infty,
\]
and hyperbolic isometries $\eta_k$ with $\eta_k(o)\to\infty$ such that,
with
\[
    W_*(z,t):=W(O_\infty z,t),
    \qquad V_*(z,t):=t^{a_0}W_*(z,t),
\]
one has
\begin{equation}\label{normalized boundary bubble}
    \bigl\|\mathcal B_k^{\mathcal T}[W]
       -\mathcal C_{\mathbb B}
          \bigl((V_*\circ\Psi^{-1})\circ\eta_k^{-1}\bigr)
    \bigr\|_{W^{m,2}_0(\mathbb B^n)}\longrightarrow0.
\end{equation}
Here $\Psi:\mathbb R^n_+\to\mathbb B^n$ is a fixed Cayley map at the
limiting boundary point, and $\mathcal C_{\mathbb B}F:=\Phi^{a_0}F$ is the
ball conformal transfer.
\end{lemma}

\begin{proof}
Use upper-half-space coordinates at the limiting boundary point and write
$c_k=\Psi(y_k)$, where $y_k=(\zeta_k,0)$ are the boundary centres.  Choose
$\Psi$ so that
\[
    \Psi^*g_{\rm Eucl}=\lambda(y)^2g_{\rm Eucl},
    \qquad \lambda(y_\infty)=1,
    \qquad y_k\longrightarrow y_\infty.
\]
The normalization of $\lambda$ is obtained by precomposing the Cayley map
with an upper-half-space dilation fixing $y_\infty$.
After choosing compatible boundary frames, conformality and preservation of
the inward normal give
\[
    D\Psi(y_k)=\lambda(y_k)Q_k,
    \qquad
    Q_k=\operatorname{diag}(O_k,1),\qquad O_k\in O(n-1),
    \qquad Q_k(\mathbb R^n_+)=\mathbb R^n_+.
\]
Pass to a subsequence for which $O_k\to O_\infty$, and set
\[
    \varepsilon_k':=\frac{\varepsilon_k}{\lambda(y_k)},
    \qquad
    h_k(z,t):=(\zeta_k+\varepsilon_k'z,\varepsilon_k't),
    \qquad
    \eta_k:=\Psi\circ h_k\circ\Psi^{-1}.
\]
The normalization of $\Psi$ implies
$\varepsilon_k'/\varepsilon_k\to1$.  The charts
\[
    \mathcal T_k(\xi)
       :=\Psi\bigl(y_k+\lambda(y_k)^{-1}Q_k^{-1}\xi\bigr)
\]
have isometric differential at the origin.  Hence
\cite[Proposition~2.1]{Mazumdar} permits replacement of the original
admissible charts by $\mathcal T_k$, up to an error tending to zero strongly
in $W^{m,2}_0(\mathbb B^n)$.

For a compactly supported smooth $W$, the bubble written with
$\mathcal T_k$ and the conformal image of $V_k(z,t):=t^{a_0}W(O_kz,t)$
have, in the variables $y=y_k+\varepsilon_k'(z,t)$, the same leading term.
Their only multiplicative discrepancy is
\[
    \left(\frac{\lambda(y_k+\varepsilon_k'(z,t))}
                 {\lambda(y_k)}\right)^{-a_0},
\]
which tends to $1$ in $C^m$ on every fixed compact subset of
$\overline{\mathbb R^n_+}$.  The scale-invariant change of variables in the
homogeneous $W^{m,2}$ norm therefore gives
\[
 \bigl\|\mathcal B_k^{\mathcal T}[W]
   -\mathcal C_{\mathbb B}
      \bigl((V_k\circ\Psi^{-1})\circ\eta_k^{-1}\bigr)
 \bigr\|_{W^{m,2}_0(\mathbb B^n)}\longrightarrow0.
\]
The cut-off tails are harmless because they escape to infinity in the
rescaled variables.  Since $O_k\to O_\infty$, the continuity of the rotation
action on $\mathcal D^{m,2}(\mathbb R^n_+)$ gives $W(O_k\,\cdot)\to W_*$.
Finally, approximate an arbitrary $W\in\mathcal D^{m,2}(\mathbb R^n_+)$ by
$C_c^\infty$ functions.  The exact conformal actions are isometries and the
chart--cut-off bubble operators are uniformly bounded by the same
proposition of Mazumdar, so the preceding estimate passes to the limit and
proves \eqref{normalized boundary bubble}.  Since
$h_k(o)=(\zeta_k,\varepsilon_k')$ leaves every compact hyperbolic set,
$\eta_k(o)\to\infty$.
\end{proof}

\begin{theorem}[Boundary profiles and hyperbolic isometries;
restatement of Theorem~\ref{intro boundary identification}]
\label{half space isometry identification}
Put $a_0=(n-2m)/2$ and endow $\mathbb R^n_+$ with the hyperbolic metric
$g_+=t^{-2}(dz^2+dt^2)$.  If $W\in\mathcal D^{m,2}(\mathbb R^n_+)$
solves \eqref{half space profile equation}, then
\begin{equation}\label{half space conformal lift}
    V(z,t):=t^{a_0}W(z,t)
\end{equation}
belongs to $\mathcal D^{m,2}(\mathbb H^n)$ and solves
\[
    P_m V=|V|^{q-2}V\quad\text{on }\mathbb H^n.
\]
The conformal lift preserves the quadratic norm, the $L^q$ norm, and
the critical energy.

For $\zeta\in\mathbb R^{n-1}$ and $\varepsilon>0$, define
\[
    h_{\zeta,\varepsilon}(z,t)
       =(\zeta+\varepsilon z,\varepsilon t).
\]
Then $h_{\zeta,\varepsilon}$ is a hyperbolic isometry and
\begin{equation}\label{dilation is hyperbolic translation}
    t^{a_0}\varepsilon^{-a_0}
    W\!\left(\frac{z-\zeta}{\varepsilon},
             \frac{t}{\varepsilon}\right)
       =V\circ h_{\zeta,\varepsilon}^{-1}(z,t).
\end{equation}
Consequently, after transporting by a fixed Cayley map between
$\mathbb R^n_+$ and $\mathbb B^n$, every boundary bubble in
Theorem~\ref{Euclidean Struwe Decomposition} is, up to a strongly
vanishing cut-off error and a fixed tangential orthogonal reparametrization
of $W$, the conformal image of
\[
    V\circ\eta_k^{-1},
    \qquad \eta_k\in\ISO(\mathbb H^n),
    \qquad \eta_k(o)\longrightarrow\infty.
\]
\end{theorem}

\begin{proof}
The first assertion is the conformal covariance of $P_m$ for
$g_+=t^{-2}g_{\mathrm{Eucl}}$; at the critical exponent it also gives
the three norm and energy identities.  Formula
\eqref{dilation is hyperbolic translation} follows directly from
\eqref{half space conformal lift}.  Moreover, for the upper-half-space
base point $o=(0,1)$,
$h_{\zeta,\varepsilon}(o)=(\zeta,\varepsilon)$ leaves every compact
hyperbolic set when $\varepsilon\to0$.

The final assertion is Lemma~\ref{Cayley chart normalization}.
After relabeling the fixed rotated profile $W_*$ supplied there as $W$, it
also shows that the scalar and orthogonal parts of the boundary-chart
normalization are absorbed, respectively, into an asymptotically equivalent
scale and a fixed tangential rotation of the profile.

If two boundary centres have different limiting boundary points, their
isometry orbits have mutual hyperbolic distance tending to infinity.  If
they have the same limiting point, the same conclusion follows from
\eqref{Mazumdar parameter orthogonality} in a common upper-half-space
chart.  This proves the last statement and the geometric separation.
\end{proof}

The following consequence for minimizing sequences will be used in
Section~\ref{sec:existence}.  We include the appeal to Ekeland's principle,
since a minimizing sequence need not be Palais--Smale to begin with.

\begin{corollary}\label{Euclidean Global Compactness}
Let $\{u_k\}\subset W^{m,2}_0(\mathbb B^n)$ satisfy
\[
    \int_{\mathbb B^n}|u_k|^q\,dx=1,
    \qquad
    \int_{\mathbb B^n}|\nabla^m u_k|^2\,dx\longrightarrow S_{m,n}.
\]
Then, after passing to a subsequence, there are $y_k\in\mathbb B^n$,
$\delta_k\to0$ with
\[
    \frac{\dist(y_k,\partial\mathbb B^n)}{\delta_k}
       \longrightarrow\infty,
\]
cut-offs $\chi_k\in C_c^\infty(\mathbb B^n)$ equal to one near $y_k$,
signs $s_k\in\{-1,1\}$, and $r_k\to0$ strongly in
$W^{m,2}_0(\mathbb B^n)$ such that
\[
    u_k=r_k+s_k\,
       \frac{\chi_kG_{\delta_k,y_k}}
            {\|\chi_kG_{\delta_k,y_k}\|_{L^q(\mathbb B^n)}},
    \qquad
    G_{\delta,y}(x)
       =\left(\frac{2\delta}{\delta^2+|x-y|^2}\right)^{\frac{n-2m}{2}}.
\]
\end{corollary}

\begin{proof}
Ekeland's variational principle gives a constrained Palais--Smale
sequence at distance $o(1)$ from $u_k$.  After rescaling by its
Lagrange multipliers, Theorem~\ref{Euclidean Struwe Decomposition}
applies.  Every nonzero profile carries at least the one-bubble energy,
so the limiting energy permits at most one nonzero object among the weak
limit and the profiles.  A nonzero weak limit at this level would attain
the sharp Sobolev constant in $W^{m,2}_0(\mathbb B^n)$; its extension by
zero would be a whole-space extremal that vanishes on an open set, which
is impossible by the equality-case classification for the higher-order
sharp Sobolev inequality \cite{Lieb,Swanson}.  Hence there is exactly
one profile.  It
cannot be a boundary profile: extension by zero would then
attain the sharp whole-space Sobolev constant while vanishing on an
open half-space, contrary to the classification of the equality cases
in the sharp Sobolev inequality.  The unique interior profile is,
up to sign, dilation and translation, the standard Sobolev extremal.
The claimed formula and strong remainder follow.
\end{proof}

\section{Global compactness on \texorpdfstring{$\mathbb H^n$}{hyperbolic space}}
\label{sec:globalcompactness}

Theorem~\ref{Global compactness} rests on two facts.  The first is that the
critical potential is compact as a map from $W^{m,2}_0(\mathbb B^n)$ to its
dual, so that it disappears once the weak limit has been removed.  The
second is Theorem~\ref{half space isometry identification}: the boundary
profiles of Mazumdar's decomposition are hyperbolic profiles transported to
infinity.

We keep the notation of Section~\ref{sec:preliminaries}.  Thus
$\Phi=2/(1-|x|^2)$, $a_0=(n-2m)/2$, $\iota$ is the ball-model
identification, and the conformal transfer
$\mathcal Cu=\Phi^{a_0}(u\circ\iota)$ of
\eqref{preliminary conformal transfer} is an isometric isomorphism from
$\mathcal D^{m,2}(\mathbb H^n)$ onto $W^{m,2}_0(\mathbb B^n)$ for the
critical quadratic forms, under which
\eqref{preliminary quadratic identity} and
\eqref{preliminary Lq identity} hold.  If
$a\in L^{n/(2m)}(\mathbb H^n)$, put
\begin{equation}\label{transformed potential}
    A(x):=(a\circ\iota)(x)\Phi(x)^{2m}.
\end{equation}
Then
\begin{equation}\label{potential conformal identities}
    \int_{\mathbb B^n}|A|^{n/(2m)}\,dx
       =\int_{\mathbb H^n}|a|^{n/(2m)}\,dV,
    \qquad
    \int_{\mathbb B^n}A(\mathcal Cu)^2\,dx
       =\int_{\mathbb H^n}au^2\,dV.
\end{equation}

\subsection{Compactness of the critical potential}

\begin{lemma}\label{critical potential compactness}
Let $A\in L^{n/(2m)}(\mathbb B^n)$.  Multiplication by $A$ defines a
compact operator
\[
    M_A:W^{m,2}_0(\mathbb B^n)\longrightarrow
        W^{-m,2}(\mathbb B^n).
\]
In particular, if $z_k\rightharpoonup0$ in
$W^{m,2}_0(\mathbb B^n)$, then
\[
    \|Az_k\|_{W^{-m,2}}\longrightarrow0,
    \qquad
    \int_{\mathbb B^n}A z_k^2\,dx\longrightarrow0.
\]
The second conclusion remains true with $A$ replaced by $|A|$.
\end{lemma}

\begin{proof}
Choose $A_j\in L^\infty(\mathbb B^n)$ such that
$A_j\to A$ in $L^{n/(2m)}$.  Rellich's theorem gives
$z_k\to0$ in $L^2(\mathbb B^n)$.  Hence, uniformly for
$\|\varphi\|_{W^{m,2}_0}\leq1$,
\[
    \left|\int_{\mathbb B^n}A_jz_k\varphi\,dx\right|
       \leq C\|A_j\|_\infty\|z_k\|_2\longrightarrow0.
\]
On the other hand, since
\[
    \frac{2m}{n}+\frac2q=1,
\]
H\"older's inequality and the critical Sobolev embedding give
\[
    \left|\int_{\mathbb B^n}(A-A_j)z_k\varphi\,dx\right|
       \leq C\|A-A_j\|_{n/(2m)}
              \|z_k\|_{W^{m,2}_0}\|\varphi\|_{W^{m,2}_0}.
\]
First let $k\to\infty$ and then $j\to\infty$.  This proves compactness.
Pairing $Az_k\to0$ in the dual with the bounded sequence $z_k$ gives
the quadratic conclusion.  The same proof applies to $|A|$.
\end{proof}

On the ball, consider
\begin{align*}
    \mathcal I_A(w)
      &:=\frac12\int_{\mathbb B^n}
          \bigl(|\nabla^m w|^2+Aw^2\bigr)\,dx
          -\frac1q\int_{\mathbb B^n}|w|^q\,dx,\\
    E_{\mathbb B^n}(w)
      &:=\frac12\int_{\mathbb B^n}|\nabla^m w|^2\,dx
          -\frac1q\int_{\mathbb B^n}|w|^q\,dx.
\end{align*}

\begin{lemma}\label{remove critical potential}
Assume $A\in L^{n/(2m)}(\mathbb B^n)$ and let $\{w_k\}$ be a bounded
Palais--Smale sequence for $\mathcal I_A$.  If
$w_k\rightharpoonup w$ in $W^{m,2}_0(\mathbb B^n)$, then $w$ is a
critical point of $\mathcal I_A$ and, with $z_k=w_k-w$,
\begin{align}
    E_{\mathbb B^n}'(z_k)&\longrightarrow0,
       \label{pure residual derivative}\\
    \mathcal I_A(w_k)&=\mathcal I_A(w)+E_{\mathbb B^n}(z_k)+o(1).
       \label{pure residual energy}
\end{align}
\end{lemma}

\begin{proof}
After taking a subsequence, $w_k\to w$ almost everywhere and locally
strongly below the critical exponent.  The standard weak convergence of
the critical Nemytskii map, together with
$Az_k\to0$ in $W^{-m,2}$ from
Lemma~\ref{critical potential compactness}, identifies $w$ as a critical
point.  The
Br\'ezis--Lieb lemma \cite{BrezisLieb1} and its derivative version yield
\[
    |w_k|^{q-2}w_k-|z_k|^{q-2}z_k-|w|^{q-2}w
       \longrightarrow0\quad\text{in }L^{q/(q-1)}.
\]
Lemma~\ref{critical potential compactness} gives
$Az_k\to0$ in $W^{-m,2}$ and $\int A z_k^2=o(1)$.  In particular,
\[
  \int_{\mathbb B^n}A(w_k^2-w^2-z_k^2)\,dx
     =2\langle Az_k,w\rangle=o(1).
\]
Subtracting the
equations for $w_k$ and $w$ proves
\eqref{pure residual derivative}; the quadratic Hilbert-space splitting
and the Br\'ezis--Lieb identity prove \eqref{pure residual energy}.
\end{proof}

\subsection{Proof of the global decomposition}

For an isometry $\eta\in\ISO(\mathbb H^n)$, we write
\[
    \eta\star V:=V\circ\eta^{-1}.
\]

\begin{proof}[Proof of Theorem~\ref{Global compactness}]
Let $\{u_k\}$ be a Palais--Smale sequence for $I$ at level $d$.
Because $a\geq0$, the identity
\[
    I(u_k)-\frac1q\langle I'(u_k),u_k\rangle
       =\frac mn\left(
          \int_{\mathbb H^n}u_kP_m u_k\,dV
          +\int_{\mathbb H^n}a u_k^2\,dV\right)
\]
shows that $\{u_k\}$ is bounded in $\mathcal D^{m,2}(\mathbb H^n)$.
Set $w_k=\mathcal Cu_k$.  By
\eqref{preliminary quadratic identity}, \eqref{preliminary Lq identity}
and \eqref{potential conformal identities}, $\{w_k\}$ is a Palais--Smale
sequence for $\mathcal I_A$ at the same level.  Passing to a subsequence, write $w_k\rightharpoonup w$ and
$z_k=w_k-w$.

Lemma~\ref{remove critical potential} shows that $w$ solves the full
ball equation and that $\{z_k\}$ is a Palais--Smale sequence for the
pure critical functional $E_{\mathbb B^n}$.  Apply
Theorem~\ref{Euclidean Struwe Decomposition}.  Since $z_k\rightharpoonup0$,
its weak profile is zero, and we obtain finitely many interior profiles
$U^\alpha\in\mathcal D^{m,2}(\mathbb R^n)$ and boundary profiles
$W^\beta\in\mathcal D^{m,2}(\mathbb R^n_+)$, together with a strong
remainder.

Let $u=\mathcal C^{-1}w$.  For an interior ball bubble
$\mathcal B_k^\alpha$, define its hyperbolic lift by
\[
    \mathcal U_k^\alpha:=\mathcal C^{-1}\mathcal B_k^\alpha.
\]
It is a cut-off concentrating bubble in ball coordinates; its scale
tends to zero and its profile $U^\alpha$ solves
\eqref{Euclidean Whole Space problem}.  For a boundary profile, let
$V^\beta$ be the conformal lift given by
Theorem~\ref{half space isometry identification}.  That theorem and the
chart-independence in Mazumdar's theorem give isometries
$\eta_k^\beta$ such that
\[
    \mathcal C^{-1}\mathcal B_k^\beta
       =\eta_k^\beta\star V^\beta+o(1)
       \quad\text{in }\mathcal D^{m,2}(\mathbb H^n),
    \qquad \eta_k^\beta(o)\to\infty.
\]
Undoing the conformal transfer in
\eqref{bounded profile decomposition}, we therefore obtain
\begin{equation}\label{final hyperbolic decomposition}
    u_k=u+\sum_{\alpha=1}^{N_{\mathrm E}}\mathcal U_k^\alpha
          +\sum_{\beta=1}^{N_{\mathrm H}}
             \eta_k^\beta\star V^\beta+r_k,
    \qquad
    r_k\to0\quad\text{in }\mathcal D^{m,2}(\mathbb H^n).
\end{equation}
The ball parameters obey
\eqref{Mazumdar parameter orthogonality}.  In particular, two boundary
profiles are separated by a hyperbolic distance tending to infinity;
in upper-half-space coordinates this follows from
\[
    \cosh d_{\mathbb H}\bigl((z_i,t_i),(z_j,t_j)\bigr)
       =\frac{|z_i-z_j|^2+t_i^2+t_j^2}{2t_it_j}.
\]

Finally, \eqref{pure residual energy},
\eqref{bounded energy splitting}, and conformal invariance give the
exact limiting identity
\[
    d=I(u)+\sum_{\alpha=1}^{N_{\mathrm E}}J(U^\alpha)
          +\sum_{\beta=1}^{N_{\mathrm H}}I_\infty(V^\beta).
\]
The corresponding quadratic and $L^q$ splittings follow from
\eqref{bounded norm splitting}--\eqref{bounded Lq splitting}.  This
proves the theorem.
\end{proof}

\section{Polar cones, energy doubling, and the compactness window}
\label{sec:moreau}

For $m\geq2$ the decomposition of a function into its positive and negative
parts leaves $W^{m,2}$, and with it the standard route to the energy bound
for sign-changing critical points.  What replaces it here is an order
property of the inverse Riesz map, which makes no reference to any
particular operator.  We isolate that property first and apply it to $P_m$
afterwards.

\subsection{Polar cones and inverse positivity}

Let $(X,\mu)$ be a measure space, let $2<q<\infty$, and let $H$ be a
real Hilbert space continuously embedded in $L^q(X,\mu)$.  Denote its
inner product by $\mathfrak a(\cdot,\cdot)$ and set
\[
    \|u\|_H^2:=\mathfrak a(u,u),
    \qquad
    S_H:=\inf_{u\in H\setminus\{0\}}
       \frac{\|u\|_H^2}{\|u\|_{L^q}^2}>0.
\]
The positive cone and its polar cone are
\[
    \mathcal K:=\{u\in H:u\geq0\ \mu\text{-a.e.}\},
    \qquad
    \mathcal K^\circ
      :=\{w\in H:\mathfrak a(w,v)\leq0
                     \text{ for every }v\in\mathcal K\}.
\]
The cone $\mathcal K$ is closed and convex.
Throughout this section, ``$u$ changes sign'' means
\[
    \mu(\{u>0\})>0
    \quad\text{and}\quad
    \mu(\{u<0\})>0.
\]

\begin{lemma}[Moreau decomposition \cite{Moreau}]
\label{Convex cone decomposition}
Every $u\in H$ admits a unique decomposition
\begin{equation}\label{abstract Moreau decomposition}
    u=p+n,
    \qquad p\in\mathcal K,
    \qquad n\in\mathcal K^\circ,
    \qquad \mathfrak a(p,n)=0.
\end{equation}
In particular,
$\|u\|_H^2=\|p\|_H^2+\|n\|_H^2$.
\end{lemma}

\begin{proof}
Let $p$ be the metric projection of $u$ onto $\mathcal K$ and put
$n=u-p$.  The variational characterization of the projection is
\[
    \mathfrak a(n,v-p)\leq0\qquad(v\in\mathcal K).
\]
Taking $v=p+tk$ with $k\in\mathcal K$, $t>0$, gives
$\mathfrak a(n,k)\leq0$, while $v=0$ and $v=2p$ give
$\mathfrak a(n,p)=0$.  Thus $n\in\mathcal K^\circ$ and the asserted
orthogonality holds.  The same variational inequalities give uniqueness.
\end{proof}

Let $\mathcal R:H\to H'$ be the Riesz map associated with
$\mathfrak a$, and let
\[
    \mathcal K'_+
       :=\{F\in H':\langle F,v\rangle\geq0
                       \text{ for every }v\in\mathcal K\}
\]
be the positive dual cone.

\begin{lemma}[Cone condition and inverse positivity]
\label{polar inverse positivity}
One has
\begin{equation}\label{Riesz polar identity}
    \mathcal R^{-1}(\mathcal K'_+)=-\mathcal K^\circ.
\end{equation}
Consequently, the following assertions are equivalent:
\begin{enumerate}[(i)]
    \item $\mathcal K^\circ\subset-\mathcal K$;
    \item $\mathcal R^{-1}(\mathcal K'_+)\subset\mathcal K$.
\end{enumerate}
If, more strongly,
\[
    \mathcal R^{-1}(\mathcal K'_+\setminus\{0\})
       \subset\{u\in\mathcal K:u>0\ \mu\text{-a.e.}\},
\]
then every nonzero $w\in\mathcal K^\circ$ is strictly negative almost
everywhere.
\end{lemma}

\begin{proof}
If $F\in\mathcal K'_+$ and $z=\mathcal R^{-1}F$, then
$\mathfrak a(z,v)\geq0$ for every $v\in\mathcal K$, so
$z\in-\mathcal K^\circ$.  The converse follows by taking $F=\mathcal Rz$.
This proves \eqref{Riesz polar identity} and the equivalence.  Under the
stronger assumption, apply it to
$-w=\mathcal R^{-1}(\mathcal R(-w))$.
\end{proof}

\subsection{Energy doubling for sign-changing critical points}

\begin{theorem}[Energy doubling for sign-changing critical points;
restatement of Theorem~\ref{intro nodal principle}]
\label{abstract nodal energy theorem}
Assume
\begin{equation}\label{abstract cone condition}
    \mathcal K^\circ\subset-\mathcal K.
\end{equation}
Consider
\[
    \mathcal E(u)
       :=\frac12\|u\|_H^2-\frac1q\int_X|u|^q\,d\mu.
\]
If $u\neq0$ is a critical point of $\mathcal E$ and changes sign, then
\begin{align}
    \|u\|_H^2=\|u\|_{L^q}^q
       &\geq2S_H^{q/(q-2)},
       \label{abstract norm doubling}\\
    \frac{\|u\|_H^2}{\|u\|_{L^q}^2}
       &\geq2^{(q-2)/q}S_H,
       \label{abstract quotient doubling}\\
    \mathcal E(u)
       &\geq2\left(\frac12-\frac1q\right)
                    S_H^{q/(q-2)}.
       \label{abstract energy doubling}
\end{align}
\end{theorem}

\begin{proof}
Let $u=p+n$ be the decomposition
\eqref{abstract Moreau decomposition}.  By
\eqref{abstract cone condition}, $p\geq0$ and $n\leq0$.  Since $u$
changes sign, neither component vanishes.  Pointwise,
\begin{equation}\label{abstract pointwise inequalities}
    |u|^{q-2}up\leq p^q,
    \qquad
    |u|^{q-2}un\leq|n|^q.
\end{equation}
Indeed, on $\{u\geq0\}$ one has $0\leq u\leq p$, while on
$\{u<0\}$ the first left-hand side is nonpositive; the second inequality
is symmetric.

Testing the equation against $p$ and $n$ and using orthogonality gives,
for $z=p,n$,
\[
    S_H\|z\|_{L^q}^2
       \leq\|z\|_H^2
       =\mathfrak a(u,z)
       \leq\|z\|_{L^q}^q.
\]
Thus $\|z\|_{L^q}^{q-2}\geq S_H$ and
$\|z\|_H^2\geq S_H^{q/(q-2)}$.  Orthogonality proves
\eqref{abstract norm doubling}.  Testing against $u$ gives
$\|u\|_H^2=\|u\|_{L^q}^q$, from which
\eqref{abstract quotient doubling} and
\eqref{abstract energy doubling} follow.
\end{proof}

\begin{remark}[In what sense the constant is optimal]
\label{abstract factor sharpness}
Within the class of spaces covered by
Theorem~\ref{abstract nodal energy theorem} the factor two cannot be
improved.  On a two-point measure space take $H=\mathbb R^2$ with its
Euclidean inner product and the coordinatewise positive cone; then
$S_H=1$, and $(1,-1)$ is a critical point for which
\eqref{abstract energy doubling} is an equality.  For a fixed connected
elliptic problem the theorem says nothing of the kind, and we make no claim
that the lower bound is attained, or even approached, by nodal solutions of
such a problem; that would call for a two-bubble construction of a quite
different nature.
\end{remark}

\subsection{Applications to the limiting equations}

\begin{proposition}\label{Positive remaining}
Let
\[
    \mathcal K=\{u\in\mathcal D^{m,2}(\mathbb H^n):u\geq0\ \text{a.e.}\}
\]
and take the inner product
\[
    \mathfrak a_0(u,v)=\int_{\mathbb H^n}uP_m v\,dV.
\]
Then $\mathcal K^\circ\subset-\mathcal K$.
\end{proposition}

\begin{proof}
Lu and Yang \cite{LuYang2} proved that the Green kernel $G_m(x,y)$ of
$P_m$ on $\mathbb H^n$ is positive.  Given
$0\leq\varphi\in C_c^\infty(\mathbb H^n)$, let
\[
    v_\varphi(x)=\int_{\mathbb H^n}G_m(x,y)\varphi(y)\,dV_y.
\]
Since $\varphi$ defines an element of the dual space and the solution of
$P_mv=\varphi$ in $\mathcal D^{m,2}(\mathbb H^n)$ is unique, the kernel
representation of \cite{LuYang2} identifies $v_\varphi$ with
$\mathcal R^{-1}\varphi$.  In particular
$v_\varphi\in\mathcal D^{m,2}(\mathbb H^n)\cap\mathcal K$.  Hence, for
$w\in\mathcal K^\circ$,
\[
    \int_{\mathbb H^n}w\varphi\,dV
       =\mathfrak a_0(w,v_\varphi)\leq0.
\]
Testing against all nonnegative compactly supported smooth functions gives
$w\leq0$ almost everywhere.
\end{proof}

\begin{corollary}\label{Energy Doubling Cor 3.3}
Every nonzero solution of \eqref{hyperbolic infinity problem} changes sign,
and satisfies
\[
    I_\infty(u)\geq\frac{2m}{n}S_{m,n}^{n/(2m)}.
\]
\end{corollary}

\begin{proof}
If a nonzero solution had one sign, change its sign if necessary so
that it is nonnegative and nontrivial.  Its conformal image is then a
nonnegative nontrivial weak solution of
\eqref{Euclidean unit ball problem}.  Critical weak solutions have the
boundary regularity required below
\cite[Proposition~7.15]{GazzolaGrunauSweers1}.  The positivity of Boggio's kernel
\cite{GazzolaGrunauSweers1} for the unperturbed ball operator
upgrades a nonnegative nontrivial solution to a strictly positive one.  The ball symmetry theorem
\cite[Theorem~7.1]{GazzolaGrunauSweers1} then makes this solution radial,
while the critical Dirichlet nonexistence theorem
\cite[Theorem~7.34]{GazzolaGrunauSweers1} excludes every nonzero radial
solution for arbitrary $m\in\mathbb N$ and $n>2m$.  Thus every nonzero
hyperbolic solution must change sign.

The energy bound now follows from
Theorem~\ref{abstract nodal energy theorem}, since
\[
    \frac12-\frac1q=\frac mn,
    \qquad
    \frac{q}{q-2}=\frac{n}{2m}.
\]
\end{proof}

\begin{proposition}[Energy of entire profiles]
\label{whole space energy alternatives}
Let $U\in\mathcal D^{m,2}(\mathbb R^n)$ be a nonzero weak solution of
\eqref{Euclidean Whole Space problem}.  Then
\[
    J(U)=\frac mnS_{m,n}^{n/(2m)}
       \quad\text{if $U$ has one sign},
    \qquad
    J(U)\geq\frac{2m}{n}S_{m,n}^{n/(2m)}
       \quad\text{if $U$ changes sign}.
\]
\end{proposition}

\begin{proof}
Put
\[
    \mathcal K_{\mathbb R^n}
       :=\{u\in\mathcal D^{m,2}(\mathbb R^n):u\geq0\ \text{a.e.}\}
\]
and endow $\mathcal D^{m,2}(\mathbb R^n)$ with the inner product
\[
    \mathfrak a_{\mathbb R^n}(u,v)
       :=\int_{\mathbb R^n}\nabla^m u\cdot\nabla^m v\,dx.
\]
The inverse Riesz map is the Riesz potential.  More precisely, for
$0\leq\varphi\in C_c^\infty(\mathbb R^n)$,
\[
    v_\varphi(x):=c_{m,n}\int_{\mathbb R^n}
       \frac{\varphi(y)}{|x-y|^{n-2m}}\,dy
\]
belongs to $\mathcal D^{m,2}(\mathbb R^n)\cap
\mathcal K_{\mathbb R^n}$ and is the Riesz solution of
$(-\Delta)^m v_\varphi=\varphi$.  Thus the argument of
Proposition~\ref{Positive remaining}, with the Riesz kernel in place of
$G_m$, gives
\[
    \mathcal K_{\mathbb R^n}^{\circ}
       \subset-\mathcal K_{\mathbb R^n}.
\]
If $U$ changes sign, Theorem~\ref{abstract nodal energy theorem} therefore
applies directly and yields the asserted two-bubble lower bound.

It remains to treat the one-sign case.  After replacing $U$ by $-U$ if
necessary, assume $U\geq0$.  Since
$U^{q-1}\in L^{q/(q-1)}(\mathbb R^n)=L^{2n/(n+2m)}(\mathbb R^n)$, the
Hardy--Littlewood--Sobolev inequality places the Riesz potential
\[
    w(x):=c_{m,n}\int_{\mathbb R^n}
       \frac{U(y)^{q-1}}{|x-y|^{n-2m}}\,dy
\]
in $\mathcal D^{m,2}(\mathbb R^n)$, and $(-\Delta)^m(U-w)=0$ weakly.
Testing this against $U-w$ gives $\|\nabla^m(U-w)\|_{L^2}=0$, so that
\begin{equation}\label{whole space integral equation}
    U(x)=c_{m,n}\int_{\mathbb R^n}
       \frac{U(y)^{q-1}}{|x-y|^{n-2m}}\,dy .
\end{equation}
As $U\geq0$ and $U\not\equiv0$, the strict positivity of the kernel in
\eqref{whole space integral equation} forces $U>0$ throughout
$\mathbb R^n$; and $\mathcal D^{m,2}(\mathbb R^n)$ embeds in
$L^{2n/(n-2m)}(\mathbb R^n)$, so $U$ is a positive regular solution of that
integral equation in the sense of Chen, Li and Ou.  By
\cite[Theorem~1]{ChenLiOu}, and independently by the method of moving
spheres \cite{LiYY}, $U$ is a standard Sobolev bubble.  At the level of the
differential equation the corresponding classification is due to Lin
\cite{Lin} for $m=2$ and to Wei and Xu \cite{WeiXu} for every integer $m$;
for a general positive solution the passage to
\eqref{whole space integral equation} rests on the super-polyharmonic
property $(-\Delta)^iU>0$, $1\leq i\leq m-1$, proved in \cite{WeiXu}, and it
is the finite-energy hypothesis that lets us reach it directly here.  Consequently
\[
    \int_{\mathbb R^n}|\nabla^mU|^2\,dx
       =\int_{\mathbb R^n}|U|^q\,dx
       =S_{m,n}^{n/(2m)},
\]
which gives the stated one-bubble energy.
\end{proof}

\subsection{The compactness window}

\begin{proposition}\label{Compactness Energy Level}
Let $a\geq0$ belong to $L^{n/(2m)}(\mathbb H^n)$.  If $\{u_k\}$ is a
Palais--Smale sequence for $I$ at a level
\[
    d\in\left(\frac mnS_{m,n}^{n/(2m)},
               \frac{2m}{n}S_{m,n}^{n/(2m)}\right),
\]
then a subsequence converges strongly in
$\mathcal D^{m,2}(\mathbb H^n)$.
\end{proposition}

\begin{proof}
Apply Theorem~\ref{Global compactness}.  Every hyperbolic profile carries
at least $\frac{2m}{n}S_{m,n}^{n/(2m)}$ by
Corollary~\ref{Energy Doubling Cor 3.3}, so none can occur.  By
Proposition~\ref{whole space energy alternatives}, an entire profile either carries
exactly the one-bubble energy or at least the two-bubble energy.

If the weak limit $u$ is nonzero, testing its equation and applying the
Sobolev inequality gives
\[
    I(u)=\frac mn\left(
       \int_{\mathbb H^n}uP_m u\,dV
       +\int_{\mathbb H^n}au^2\,dV\right)
       \geq\frac mnS_{m,n}^{n/(2m)}.
\]
Thus a nonzero weak limit cannot coexist with an entire profile below the
two-bubble level, and two entire profiles are excluded for the same reason.
A single sign-changing entire profile is excluded as well, since by
Proposition~\ref{whole space energy alternatives} its energy already
exceeds $d$.
The only remaining noncompact possibility is one one-sign entire profile
and zero weak limit, which has level exactly
$\frac mnS_{m,n}^{n/(2m)}$, contrary to the strict lower bound on $d$.
Hence no profile occurs, and the remainder in
Theorem~\ref{Global compactness} converges strongly.
\end{proof}

It remains to pass from constrained to unconstrained Palais--Smale
sequences.  The following is the form in which the compactness window will
be used in Section~\ref{sec:existence}.

\begin{corollary}[Constrained Palais--Smale window]
\label{constrained PS bridge}
Let $A\geq0$ belong to $L^{n/(2m)}(\mathbb B^n)$ and set
\[
    F_A(v)=\int_{\mathbb B^n}(|\nabla^m v|^2+Av^2)\,dx,
    \qquad
    \mathcal M=\left\{v\in W^{m,2}_0(\mathbb B^n):
                 \int_{\mathbb B^n}|v|^q\,dx=1\right\}.
\]
Then $F_A|_{\mathcal M}$ satisfies the Palais--Smale condition at every
level
\[
    c\in(S_{m,n},2^{2m/n}S_{m,n}).
\]
\end{corollary}

\begin{proof}
Let $v_k\in\mathcal M$ be a constrained Palais--Smale sequence with
$F_A(v_k)\to c$.  There are multipliers $\mu_k$ such that
\[
    (-\Delta)^mv_k+Av_k
       -\mu_k|v_k|^{q-2}v_k\longrightarrow0
       \quad\text{in }W^{-m,2},
    \qquad
    \mu_k=c+o(1),
\]
where the second relation follows by testing against $v_k$.  Since
$c>S_{m,n}>0$, one has $\mu_k>0$ for large $k$, and
$\mu_k^{1/(q-2)}\to c^{1/(q-2)}>0$.  Put $w_k=\mu_k^{1/(q-2)}v_k$.  Then $w_k$ is Palais--Smale for the
unconstrained ball functional and
\[
    \mathcal I_A(w_k)\longrightarrow
       \left(\frac12-\frac1q\right)c^{q/(q-2)}
       =\frac mn c^{n/(2m)}.
\]
This level lies strictly between the one- and two-bubble energy levels.
Conformal transfer to $\mathbb H^n$ and
Proposition~\ref{Compactness Energy Level} give strong convergence of
$w_k$, and hence of $v_k$.
\end{proof}

\section{A two-level minimax criterion and concentrated potentials}
\label{sec:existence}

\subsection{A barycenter--concentration criterion}

Let $\mathcal M$ be a complete $C^1$ Finsler manifold, let
$F\in C^1(\mathcal M)$, and let
\[
    \Psi=(\beta,\gamma):\mathcal M\longrightarrow
       \mathbb R^N\times\mathbb R
\]
be continuous.  Fix $\gamma_*>0$, put $z_*=(0,\gamma_*)$, and define
\[
    \mathcal A_*:=\Psi^{-1}\bigl(\{0\}\times[\gamma_*,\infty)\bigr),
    \qquad
    \mathcal Z_*:=\Psi^{-1}(z_*).
\]
Let $K$ be the closure of a bounded open subset of $\mathbb R^{N+1}$
which is homeomorphic to a ball, and let $\Pi:K\to\mathcal M$ be
continuous.  Suppose
\begin{equation}\label{abstract nonzero degree}
    z_*\notin(\Psi\circ\Pi)(\partial K),
    \qquad
    \deg(\Psi\circ\Pi,\operatorname{int}K,z_*)\neq0.
\end{equation}
Set
\begin{equation}\label{abstract four levels}
    a:=\inf_{\mathcal A_*}F,
    \quad b:=\max_{\partial K}F\circ\Pi,
    \quad d:=\inf_{\mathcal Z_*}F,
    \quad e:=\max_KF\circ\Pi.
\end{equation}

\begin{theorem}[Two-level barycenter--concentration criterion]
\label{abstract two level criterion}
Assume that $F$ satisfies the Palais--Smale condition at every level in
$(S_0,\kappa)$, where $S_0<\kappa$.
\begin{enumerate}[(i)]
\item If
    \begin{equation}\label{abstract first chain}
        S_0<a\leq b<\min\{d,\kappa\},
    \end{equation}
then $F$ has a critical point $v_1$ with
    \[
        a\leq F(v_1)\leq b.
    \]
\item If, in addition,
    \begin{equation}\label{abstract second chain}
        b<d\leq e<\kappa,
    \end{equation}
then $F$ has a second critical point $v_2$ with
    \[
        d\leq F(v_2)\leq e.
    \]
In particular, $F(v_1)<F(v_2)$.
\end{enumerate}
\end{theorem}

\begin{proof}
Suppose first that there is no critical value in $[a,b]$.  The
deformation lemma, the Palais--Smale condition, and
\eqref{abstract first chain} give a deformation of $F^b$ into
$F^{a-\delta}$ for some $\delta>0$.  Apply it to $\Pi|_{\partial K}$.
During the deformation the $\Psi$-image cannot meet $z_*$, because the
energy remains at most $b<d$.  At the end it cannot meet the ray
$\{0\}\times[\gamma_*,\infty)$, by the definition of $a$.  The complement
of this ray is contractible inside $\mathbb R^{N+1}\setminus\{z_*\}$:
after translating $z_*$ to the origin and rotating the ray to the last
coordinate axis, polar coordinates identify the complement with
$(0,\infty)\times(S^N\setminus\{e_{N+1}\})$.  Thus
$\Psi\circ\Pi|_{\partial K}$ would be null-homotopic in
$\mathbb R^{N+1}\setminus\{z_*\}$, contradicting
\eqref{abstract nonzero degree}.  This proves (i).

For the second point, let
\[
    \Gamma:=\{h\in C(K,\mathcal M):
                    h|_{\partial K}=\Pi|_{\partial K}\},
    \qquad
    c_*:=\inf_{h\in\Gamma}\max_{\xi\in K}F(h(\xi)).
\]
For every $h\in\Gamma$, homotopy invariance of the Brouwer degree gives
\[
    \deg(\Psi\circ h,\operatorname{int}K,z_*)
       =\deg(\Psi\circ\Pi,\operatorname{int}K,z_*)\neq0.
\]
Hence $h(K)$ meets $\mathcal Z_*$, and therefore
\[
    d\leq c_*\leq e.
\]
Since the fixed boundary has energy at most $b<d$, the relative
deformation lemma applies without moving it.  The Palais--Smale condition
at $c_*\in[d,e]\subset(S_0,\kappa)$ yields a critical point $v_2$ at
level $c_*$.  Its level is strictly above that of $v_1$.
\end{proof}

We turn to the application.  By the conformal equivalence recalled in
Section~\ref{GJMS eigenvalues} it is enough to treat the Euclidean
problem \eqref{Euclidean equation} on $\Omega:=\mathbb B^n$, and we work
with the constrained functional
\[
    F_A(u) = \int_{\Omega}\left[|\nabla^m u|^2 + A(x)\,u^2\right] dx
\]
on the manifold
\[
    \mathcal N(\Omega) = \left\{u \in W_0^{m,2}(\Omega) : \int_{\Omega} |u|^q\,dx = 1\right\},
\]
where, as in Section~\ref{sec:globalcompactness},
$A:=\Phi^{2m}\,(a\circ\iota)$ is the
transferred potential.  The set $\mathcal N(\Omega)$ is a complete $C^1$
submanifold of $W^{m,2}_0(\Omega)$.  The condition on $a$ in
Theorem~\ref{Main Theorem} translates to
\[
    A(x) = \bar\alpha(x) + \lambda^{2m}\alpha[\lambda(x - x_0)],
\]
with $\bar\alpha \in L^{n/(2m)}(\Omega)$ and
$\alpha \in L^{n/(2m)}(\mathbb R^n)$ nonnegative and
$\alpha\not\equiv0$; the corresponding functional is then the
$f_\lambda$ of Theorem~\ref{Main Theorem}.

For the whole-space constructions below, we use the explicitly
higher-order normalized manifold
\[
    \mathcal N_m(\mathbb R^n):=\left\{u\in\mathcal D^{m,2}(\mathbb R^n):
             \int_{\mathbb R^n}|u|^q\,dx=1\right\}.
\]

\subsection{Barycenter and concentration maps}

On the manifold $\mathcal N_m(\mathbb R^n)$ we define
\begin{align}\label{def beta}
    \beta(u) &= \int_{\mathbb R^n} \frac{x}{1+|x|}|u(x)|^q\,dx, \\
    \label{def gamma}
    \gamma(u) &= \int_{\mathbb R^n} \left|\frac{x}{1+|x|} - \beta(u)\right||u(x)|^q\,dx.
\end{align}
Both are continuous.  The first, $\beta:\mathcal N_m(\mathbb R^n)\to\mathbb R^n$,
is a barycenter attached to the density $|u|^q$; the second,
$\gamma:\mathcal N_m(\mathbb R^n)\to\mathbb R^+$, measures how far that density
spreads away from $\beta(u)$.

\begin{remark}\label{scaling remark}
A function $u$ solves $(-\Delta)^m u + \alpha(x)u = |u|^{q-2}u$ in
$\Omega$ if and only if
$u_\lambda(x) := \lambda^{(n-2m)/2}u[\lambda(x-x_0)]$ solves
    \[
        (-\Delta)^m u_\lambda + \lambda^{2m}\alpha[\lambda(x-x_0)]u_\lambda = |u_\lambda|^{q-2}u_\lambda
    \]
in $\Omega_\lambda = x_0 + \frac{1}{\lambda}\Omega$.  Note that
$\|u_\lambda\|_{L^q(\Omega_\lambda)} = \|u\|_{L^q(\Omega)}$.  Moreover,
setting $\alpha_\lambda(x) = \lambda^{2m}\alpha[\lambda(x-x_0)]$, for
every $\epsilon > 0$:
    \[
        \int_{B(x_0,\epsilon)} \alpha_\lambda^{n/(2m)}(x)\,dx = \int_{B(0,\lambda\epsilon)} \alpha^{n/(2m)}(x)\,dx,
    \]
so that
$\lim_{\lambda \to \infty} \int_\Omega \alpha_\lambda^{n/(2m)}\,dx = \int_{\mathbb R^n} \alpha^{n/(2m)}\,dx$.
\end{remark}

For a nonnegative function $\alpha \in L^{n/(2m)}(\mathbb R^n)$, define
the critical energy level
\begin{equation}\label{def c alpha}
    c(\alpha) := \inf\left\{\int_{\mathbb R^n}\left[|\nabla^m u|^2 + \alpha(x)u^2\right]dx : u \in \mathcal N_m(\mathbb R^n),\; \beta(u) = 0,\; \gamma(u) = \tfrac13\right\}.
\end{equation}

\subsection{\texorpdfstring{The strict inequality $c(\alpha)>S_{m,n}$}
{The strict inequality c(alpha) > S(m,n)}}

\begin{lemma}\label{strict inequality c alpha}
Let $\alpha \geq 0$, $\alpha \in L^{n/(2m)}(\mathbb R^n)$, with
$\|\alpha\|_{L^{n/(2m)}} \neq 0$.  Then $c(\alpha) > S_{m,n}$.
\end{lemma}

\begin{proof}
Clearly $c(\alpha) \geq S_{m,n}$.  We show equality cannot hold.  If it
did, there would exist $(u_i) \subset \mathcal N_m(\mathbb R^n)$ with
    \begin{equation}\label{approaching S}
        \lim_{i \to \infty} \int_{\mathbb R^n}\left[|\nabla^m u_i|^2 + \alpha(x)u_i^2\right]dx = S_{m,n},
    \end{equation}
    \begin{equation}\label{barycenter constraint}
        \beta(u_i) = 0 \quad \text{and} \quad \gamma(u_i) = \tfrac13 \quad \forall\, i \in \mathbb N.
    \end{equation}
Since $\alpha \geq 0$,
    \begin{equation}\label{gradient approaching S}
        \lim_{i \to \infty} \int_{\mathbb R^n} |\nabla^m u_i|^2\,dx = S_{m,n}.
    \end{equation}
As $(u_i)$ is a minimizing sequence for the Sobolev quotient on
$\mathbb R^n$ with $\|u_i\|_{L^q} = 1$, the sharp
concentration--compactness theorem and the classification of the
higher-order Sobolev extremals \cite{Lions1,Lions2,Lieb,Swanson} give
sequences $(y_i) \subset \mathbb R^n$,
$(\sigma_i) \subset \mathbb R^+$, signs $s_i\in\{-1,1\}$, and
$(w_i) \subset \mathcal D^{m,2}(\mathbb R^n)$ such that
    \[
        u_i = w_i + s_i\frac{G_{\sigma_i, y_i}}{\|G_{\sigma_i, y_i}\|_{L^q(\mathbb R^n)}},
    \]
where $w_i \to 0$ strongly in $L^q(\mathbb R^n)$.

We claim that $(y_i)$ and $(\sigma_i)$ are bounded.

\emph{Boundedness of $(y_i)$.} Suppose $|y_i| \to +\infty$ along a
subsequence.  Set
$\Sigma_i = \{x \in \mathbb R^n : (x - y_i) \cdot y_i \geq 0\}$.  Since
$|x|/(1+|x|) \geq |y_i|/(1+|y_i|)$ for all $x \in \Sigma_i$ and
$\lim_{i \to \infty} \int_{\Sigma_i} |u_i|^q\,dx = \frac12$, we get
    \[
        \gamma(u_i) \geq \int_{\Sigma_i} \frac{|x|}{1+|x|}|u_i|^q\,dx \geq \frac{|y_i|}{1+|y_i|} \cdot \frac12,
    \]
so $\liminf_{i \to \infty} \gamma(u_i) \geq \frac12$, contradicting
\eqref{barycenter constraint}.

\emph{Boundedness of $(\sigma_i)$.} Suppose $\sigma_i \to +\infty$.  Then
$\sup_{x} |G_{\sigma_i,y_i}(x)| / \|G_{\sigma_i,y_i}\|_{L^q} \to 0$, so
$\int_{B(0,r)} |u_i|^q\,dx \to 0$ for all $r > 0$.  Hence
    \[
        \gamma(u_i) \geq \int_{\mathbb R^n \setminus B(0,r)} \frac{|x|}{1+|x|}|u_i|^q\,dx \geq \frac{r}{1+r}\int_{\mathbb R^n \setminus B(0,r)} |u_i|^q\,dx \xrightarrow{i \to \infty} \frac{r}{1+r},
    \]
which gives $\lim \gamma(u_i) = 1$, contradicting
\eqref{barycenter constraint}.

Thus, up to a subsequence, $y_i \to \bar{y} \in \mathbb R^n$ and
$\sigma_i \to \bar{\sigma} \geq 0$.  We must have $\bar{\sigma} > 0$: if
$\bar{\sigma} = 0$, then $\beta(u_i) \to \bar{y}/(1+|\bar{y}|)$, so
$\bar{y} = 0$ (since $\beta(u_i) = 0$), and then
$\gamma(u_i) \to |\bar{y}|/(1+|\bar{y}|) = 0$, contradicting
$\gamma(u_i) = \frac13$.

Passing to a further subsequence, $s_i$ is constant.  Therefore $u_i$
converges, up to sign, to
$G_{\bar{\sigma},\bar{y}} / \|G_{\bar{\sigma},\bar{y}}\|_{L^q}$ strongly
in $L^q(\mathbb R^n)$ with $\bar{\sigma} > 0$.  Since
$G_{\bar{\sigma},\bar{y}} > 0$ everywhere and
$\|\alpha\|_{L^{n/(2m)}} \neq 0$,
    \begin{equation}\label{alpha positive contribution}
        \int_{\mathbb R^n} \alpha(x)\frac{G_{\bar{\sigma},\bar{y}}^2}{\|G_{\bar{\sigma},\bar{y}}\|_{L^q}^2}\,dx > 0.
    \end{equation}
Strong $L^q$ convergence gives strong $L^{q/2}$ convergence of the
squares.  Since $\alpha\in L^{q/(q-2)}$, H\"older's inequality therefore
gives convergence of the potential integrals.  Combining this fact with
\eqref{gradient approaching S} and \eqref{alpha positive contribution}
yields
    \[
        \lim_{i \to \infty} \int_{\mathbb R^n}\left[|\nabla^m u_i|^2 + \alpha(x)u_i^2\right]dx = S_{m,n} + \int_{\mathbb R^n} \alpha(x)\frac{G_{\bar{\sigma},\bar{y}}^2}{\|G_{\bar{\sigma},\bar{y}}\|_{L^q}^2}\,dx > S_{m,n},
    \]
contradicting \eqref{approaching S}.
\end{proof}

\subsection{Concentration behavior of test functions}

For $\varphi \in W_0^{m,2}(B(0,1))$ and $\sigma > 0$,
$y \in \mathbb R^n$, define
\begin{equation}\label{def T sigma y}
    T_{\sigma,y}(u)(x) = \frac{u_{\sigma,y}}{\|u_{\sigma,y}\|_{L^q}}, \quad \text{where } u_{\sigma,y}(x) = u\!\left(\frac{x-y}{\sigma}\right).
\end{equation}
Throughout, we fix $\varphi \in W_0^{m,2}(B(0,1))$ satisfying
\begin{equation}\label{test function properties}
    \left\{\begin{array}{l}
        \varphi \in C^\infty(\overline{B(0,1)})\cap W^{m,2}_0(B(0,1)),
        \quad \varphi(x) > 0 \quad \forall\, x \in B(0,1), \\
        \varphi \text{ is radially symmetric and strictly decreasing in }|x|, \\
        \displaystyle\int_{B(0,1)} |\varphi|^q\,dx = 1, \quad S_{m,n} < \int_{B(0,1)} |\nabla^m \varphi|^2\,dx < S_{m,n} + \epsilon,
    \end{array}\right.
\end{equation}
Here $\epsilon>0$ is chosen in the proof of Theorem~\ref{Main Theorem}
so that
\[
    S_{m,n}+\epsilon
       <\min\{c(\alpha),\,2^{2m/n}S_{m,n}\}.
\]
Such $\varphi$ exists by the approximation properties of $S_{m,n}$.

\begin{lemma}\label{gamma limits}
The following relations hold:
    \begin{equation}\label{gamma properties}
        \left\{\begin{array}{ll}
            \mathrm{a)} & \displaystyle\lim_{\sigma \to 0} \sup\left\{\gamma \circ T_{\sigma,y}(\varphi) : y \in \mathbb R^n\right\} = 0, \\[6pt]
            \mathrm{b)} & (\beta \circ T_{\sigma,y}(\varphi) \cdot y) > 0 \quad \forall\, y \neq 0,\; \forall\, \sigma > 0, \\[6pt]
            \mathrm{c)} & \displaystyle\lim_{\sigma \to +\infty} \inf\left\{\gamma \circ T_{\sigma,y}(\varphi) : y \in \mathbb R^n,\, |y| \leq r\right\} = 1 \quad \forall\, r \geq 0.
        \end{array}\right.
\end{equation}
\end{lemma}

\begin{proof}
\emph{Proof of} a).  Suppose by contradiction that there exist $(y_i)$ in
$\mathbb R^n$ and $(\sigma_i)$ with $\sigma_i \to 0$ such that
$\lim_{i \to \infty} \gamma \circ T_{\sigma_i,y_i}(\varphi) > 0$.  Using
the Lipschitz estimate
    \begin{equation}\label{Lipschitz estimate}
        \left|\frac{x}{1+|x|} - \frac{y}{1+|y|}\right| \leq |x - y| \quad \forall\, x, y \in \mathbb R^n,
    \end{equation}
and the fact that $T_{\sigma_i,y_i}(\varphi)$ is supported in
$B(y_i, \sigma_i)$, we obtain
    \[
        \gamma \circ T_{\sigma_i,y_i}(\varphi) \leq \int_{B(y_i,\sigma_i)} |x - y_i|\, T_{\sigma_i,y_i}^q(\varphi)\,dx + \left|\frac{y_i}{1+|y_i|} - \beta \circ T_{\sigma_i,y_i}(\varphi)\right| \leq 2\sigma_i \to 0,
    \]
a contradiction.

\emph{Proof of} b).  By the radial symmetry and monotonicity of
$\varphi$, for $y \neq 0$,
    \[
        T_{\sigma,y}(\varphi)(x) > T_{\sigma,y}(\varphi)(-x) \quad \text{whenever } x \cdot y > 0 \text{ and } x \in B(y,\sigma).
    \]
Therefore
$(\beta \circ T_{\sigma,y}(\varphi) \cdot y) = \int_{\mathbb R^n} \frac{x \cdot y}{1+|x|}T_{\sigma,y}^q(\varphi)\,dx > 0$.

\emph{Proof of} c).  We show both
    \begin{equation}\label{limsup leq 1}
        \limsup_{\sigma \to +\infty} \inf\{\gamma \circ T_{\sigma,y}(\varphi) : |y| \leq r\} \leq 1
    \end{equation}
and
    \begin{equation}\label{liminf geq 1}
        \liminf_{\sigma \to +\infty} \inf\{\gamma \circ T_{\sigma,y}(\varphi) : |y| \leq r\} \geq 1.
    \end{equation}

For \eqref{limsup leq 1}, choose $y=0$.  Radial symmetry gives
$\beta(T_{\sigma,0}\varphi)=0$, and hence
$\gamma(T_{\sigma,0}\varphi)\leq1$.

For \eqref{liminf geq 1}, let $\sigma_i\to\infty$ and $|y_i|\leq r$.
Since the function $x/(1+|x|)$ is bounded by one, translation and
scaling give
    \[
      |\beta(T_{\sigma_i,y_i}\varphi)|
       \leq\int_{\mathbb R^n}
          \bigl|T_{1,y_i/\sigma_i}^q(\varphi)
                    -T_{1,0}^q(\varphi)\bigr|\,dx
       \longrightarrow0
    \]
uniformly for $|y_i|\leq r$.  Moreover, for fixed $\rho>0$, boundedness
of $\varphi$ gives
    \[
       \sup_{|y|\leq r}
       \int_{B(0,\rho)}T_{\sigma_i,y}^q(\varphi)\,dx
          \leq C_\varphi |B(0,\rho)|\,\sigma_i^{-n}
          \longrightarrow0.
    \]
Therefore
    \[
      \gamma(T_{\sigma_i,y_i}\varphi)
       \geq\frac{\rho}{1+\rho}
          \int_{\mathbb R^n\setminus B(0,\rho)}
                  T_{\sigma_i,y_i}^q(\varphi)\,dx
          -|\beta(T_{\sigma_i,y_i}\varphi)|.
    \]
First let $i\to\infty$ and then $\rho\to\infty$ to obtain
\eqref{liminf geq 1}.
\end{proof}

\begin{lemma}\label{potential limits}
Let $\alpha \geq 0$, $\alpha \in L^{n/(2m)}(\mathbb R^n)$.  Then:
    \begin{equation}\label{alpha vanishing}
        \left\{\begin{array}{ll}
            \mathrm{a)} & \displaystyle\lim_{\sigma \to 0} \sup\left\{\int_{\mathbb R^n} \alpha(x)T_{\sigma,y}^2(\varphi)\,dx : y \in \mathbb R^n\right\} = 0, \\[6pt]
            \mathrm{b)} & \displaystyle\lim_{\sigma \to +\infty} \sup\left\{\int_{\mathbb R^n} \alpha(x)T_{\sigma,y}^2(\varphi)\,dx : y \in \mathbb R^n\right\} = 0, \\[6pt]
            \mathrm{c)} & \displaystyle\lim_{r \to +\infty} \sup\left\{\int_{\mathbb R^n} \alpha(x)T_{\sigma,y}^2(\varphi)\,dx : \sigma > 0,\, |y| = r\right\} = 0.
        \end{array}\right.
\end{equation}
\end{lemma}

\begin{proof}
In each case, we use H\"older's inequality: since
$T_{\sigma,y}(\varphi)$ is supported in $B(y,\sigma)$ with
$\|T_{\sigma,y}(\varphi)\|_{L^q} = 1$,
    \begin{equation}\label{holder for alpha}
        \int_{\mathbb R^n} \alpha\,T_{\sigma,y}^2(\varphi)\,dx = \int_{B(y,\sigma)} \alpha\,T_{\sigma,y}^2(\varphi)\,dx \leq \left(\int_{B(y,\sigma)} \alpha^{n/(2m)}\,dx\right)^{2m/n}.
    \end{equation}

\emph{Proof of} a).  If $\sigma_i \to 0$, then $|B(y_i,\sigma_i)| \to 0$,
so $\int_{B(y_i,\sigma_i)} \alpha^{n/(2m)}\,dx \to 0$ uniformly in $y_i$
(since $\alpha^{n/(2m)}$ is integrable).

\emph{Proof of} b).  If $\sigma_i \to +\infty$, we split
    \[
        \int_{\mathbb R^n} \alpha\,T_{\sigma_i,y_i}^2\,dx \leq \left(\int_{B(0,\rho)} \alpha^{n/(2m)}\right)^{2m/n}\left(\int_{B(0,\rho)} T_{\sigma_i,y_i}^q\right)^{2/q} + \left(\int_{\mathbb R^n \setminus B(0,\rho)} \alpha^{n/(2m)}\right)^{2m/n}.
    \]
Since $\varphi$ is bounded and $L^q$-normalized,
    \[
       \sup_{y\in\mathbb R^n}
       \int_{B(0,\rho)}T_{\sigma,y}^q(\varphi)\,dx
          \leq C_\varphi |B(0,\rho)|\,\sigma^{-n}.
    \]
Thus, for fixed $\rho$, the first term tends to zero uniformly in $y_i$
as $\sigma_i\to\infty$, while the second tends to zero as
$\rho\to\infty$.

\emph{Proof of} c).  If the assertion failed, there would be
$|y_i|\to\infty$, $\sigma_i>0$, and a fixed positive lower bound for the
potential integrals.  The uniform assertions in a) and b) then force,
after passing to a subsequence,
$0<\sigma_-\leq\sigma_i\leq\sigma_+<\infty$.  Consequently,
    \[
       \int_{B(y_i,\sigma_i)}\alpha^{n/(2m)}\,dx
          \leq\int_{B(y_i,\sigma_+)}
                    \alpha^{n/(2m)}\,dx\longrightarrow0
    \]
because $\alpha^{n/(2m)}\in L^1(\mathbb R^n)$.  This contradicts
\eqref{holder for alpha}.
\end{proof}

\subsection{Topological linking construction}

\begin{corollary}\label{test family degree}
Suppose $\|\alpha\|_{L^{n/(2m)}} \neq 0$.  Then there exist $r > 0$ and
$0 < \sigma_1 < \frac13 < \sigma_2$ such that, with
$K = K(\sigma_1,\sigma_2,r) := \{(y,\sigma) : |y| \leq r,\, \sigma_1 \leq \sigma \leq \sigma_2\}$,
    \begin{equation}\label{energy on boundary K}
        \sup\left\{\int_{\mathbb R^n}\left[|\nabla^m T_{\sigma,y}(\varphi)|^2 + \alpha(x)T_{\sigma,y}^2(\varphi)\right]dx : (y,\sigma) \in \partial K\right\} < S_{m,n} + \epsilon < c(\alpha).
    \end{equation}
Moreover, the map $\Theta:K\to \mathbb R^n \times \mathbb R$ defined
by
    \[
        \Theta(y,\sigma) = \left(\beta \circ T_{\sigma,y}(\varphi),\; \gamma \circ T_{\sigma,y}(\varphi)\right)
    \]
satisfies
    \begin{equation}\label{degree of test family}
        (0,\tfrac13)\notin\Theta(\partial K),
        \qquad
        \deg\bigl(\Theta,\operatorname{int}K,(0,\tfrac13)\bigr)=1.
    \end{equation}
\end{corollary}

\begin{proof}
By \eqref{gamma properties}\,a) and \eqref{alpha vanishing}\,a), choose
$\sigma_1 \in (0,\frac13)$ such that
$\gamma \circ T_{\sigma_1,y}(\varphi) < \frac13$ for all $y$, and
the energy bound \eqref{energy on boundary K} holds when
$\sigma = \sigma_1$.  By \eqref{alpha vanishing}\,c), choose $r > 0$ such
that \eqref{energy on boundary K} holds when $|y| = r$.  By
\eqref{gamma properties}\,c) and \eqref{alpha vanishing}\,b), choose
$\sigma_2 > \frac13$ such that
$\gamma \circ T_{\sigma_2,y}(\varphi) > \frac13$ for $|y| \leq r$,
and \eqref{energy on boundary K} holds when $\sigma = \sigma_2$.

For the homotopy, define
$\vartheta: \partial K \times [0,1] \to \mathbb R^n \times \mathbb R \setminus \{(0,\frac13)\}$
by
    \[
        \vartheta(y,\sigma,t) = (1-t)(y,\sigma) + t\,\Theta(y,\sigma).
    \]
This avoids $(0,\frac13)$: when $\sigma = \sigma_1$ and
$|y| \leq r$, the $\sigma$-component
$(1-t)\sigma_1 + t\,\gamma \circ T_{\sigma_1,y}(\varphi) < \frac13$;
when $\sigma = \sigma_2$ and $|y| \leq r$, the $\sigma$-component
$> \frac13$; when $|y| = r$, by \eqref{gamma properties}\,b), the
$y$-component satisfies
$[(1-t)y + t\,\beta \circ T_{\sigma,y}(\varphi)] \cdot y > 0$.  Thus
$\Theta|_{\partial K}$ is homotopic to the identity map on $\partial K$
in the complement of $(0,\frac13)$.  Since
$(0,\frac13)\in\operatorname{int}K$, homotopy invariance of the Brouwer
degree gives \eqref{degree of test family}.
\end{proof}

\subsection{Transfer to the bounded domain}

For $\lambda > 0$, set
$\beta_\lambda = \beta \circ T_{\lambda,-\lambda x_0}$,
$\gamma_\lambda = \gamma \circ T_{\lambda,-\lambda x_0}$, and define
$f_\lambda: W_0^{m,2}(\Omega) \to \mathbb R$ by
\[
    f_\lambda(u) = \int_\Omega \left\{|\nabla^m u|^2 + \left[\bar\alpha(x) + \lambda^{2m}\alpha(\lambda(x-x_0))\right]u^2\right\}dx.
\]

\begin{lemma}\label{levels on the constraint}
For every $\lambda > 0$:
    \begin{equation}\label{inf on constraint}
        \left\{\begin{array}{ll}
            \mathrm{a)} & \inf\left\{f_\lambda(u) : u \in \mathcal N(\Omega),\, \beta_\lambda(u) = 0,\, \gamma_\lambda(u) = \tfrac13\right\} \geq c(\alpha) > S_{m,n}, \\[4pt]
            \mathrm{b)} & \inf\left\{f_\lambda(u) : u \in \mathcal N(\Omega),\, \beta_\lambda(u) = 0,\, \gamma_\lambda(u) \geq \tfrac13\right\} > S_{m,n}.
        \end{array}\right.
\end{equation}
\end{lemma}

\begin{proof}
\emph{Part} a).  For $u \in W_0^{m,2}(\Omega)$ (extended by zero), set
$u_\lambda = T_{\lambda,-\lambda x_0}(u)$.  Then
$\beta_\lambda(u) = \beta(u_\lambda)$ and
$\gamma_\lambda(u) = \gamma(u_\lambda)$.  Since $\bar\alpha \geq 0$,
    \[
        f_\lambda(u) \geq \int_{\mathbb R^n}\left[|\nabla^m u|^2 + \lambda^{2m}\alpha(\lambda(x-x_0))u^2\right]dx = \int_{\mathbb R^n}\left[|\nabla^m u_\lambda|^2 + \alpha(x)u_\lambda^2\right]dx \geq c(\alpha).
    \]

\emph{Part} b).  Suppose by contradiction that there exists
$(u_i) \subset \mathcal N(\Omega)$ with $\beta_\lambda(u_i) = 0$,
$\gamma_\lambda(u_i) \geq \frac13$, and
$f_\lambda(u_i) \to S_{m,n}$.  Since $\bar\alpha, \alpha \geq 0$, we
get $\int_\Omega |\nabla^m u_i|^2\,dx \to S_{m,n}$.

By Corollary~\ref{Euclidean Global Compactness}, there exist
$\delta_i\to0$, $x_i\in\mathbb B^n$, cut-offs $\chi_i$, signs
$s_i\in\{-1,1\}$, and $w_i\to0$ in $L^q$ such that
    \[
        u_i = w_i+s_i\frac{\chi_iG_{\delta_i,x_i}}
                    {\|\chi_iG_{\delta_i,x_i}\|_{L^q(\mathbb B^n)}}.
    \]
Setting $v_i = T_{\lambda,-\lambda x_0}(u_i)$ and
$z_i:=\lambda(x_i-x_0)$, we have $\beta(v_i)=0$ and
$\gamma(v_i)\geq\frac13$.  The sign disappears from the density, while
strong $L^q$ convergence of the remainder implies that the $L^1$
difference of the corresponding $q$-densities tends to zero.  Since the
rescaled bubble has centre $z_i$ and scale $\lambda\delta_i\to0$, it
follows that
    \begin{equation}\label{concentration of vi}
        \int_{\mathbb R^n\setminus B(z_i,\rho)}|v_i|^q\,dx
             \longrightarrow0
        \qquad\text{for every fixed }\rho>0.
    \end{equation}
Put $g(x):=x/(1+|x|)$.  From $\beta(v_i)=0$, the Lipschitz estimate
\eqref{Lipschitz estimate}, and \eqref{concentration of vi},
    \[
       |g(z_i)|=|g(z_i)-\beta(v_i)|
       \leq \rho
          +2\int_{\mathbb R^n\setminus B(z_i,\rho)}|v_i|^q\,dx.
    \]
Taking $i\to\infty$ and then $\rho\downarrow0$ gives $g(z_i)\to0$, hence
$z_i\to0$.  Finally,
    \[
      \gamma(v_i)
       \leq\int_{\mathbb R^n}|g(x)-g(z_i)|\,|v_i|^q\,dx
             +|g(z_i)-\beta(v_i)|\longrightarrow0,
    \]
again by first restricting to $B(z_i,\rho)$ and then letting
$\rho\downarrow0$.  This contradicts $\gamma(v_i)\geq\frac13$.
\end{proof}

\begin{lemma}\label{rescaled test family bounds}
There exists $\bar{\lambda} > 0$ such that for every
$\lambda > \bar{\lambda}$:
    \begin{equation}\label{support and energy}
        \left\{\begin{array}{ll}
            \mathrm{a)} & \operatorname{supp}\left(T_{1/\lambda,x_0} \circ T_{\sigma,y}(\varphi)\right) \subset \Omega \quad \forall\, (y,\sigma) \in K, \\[4pt]
            \mathrm{b)} & \displaystyle\sup\left\{f_\lambda \circ T_{1/\lambda,x_0} \circ T_{\sigma,y}(\varphi) : (y,\sigma) \in \partial K\right\} < S_{m,n} + \epsilon < c(\alpha).
        \end{array}\right.
\end{equation}
\end{lemma}

\begin{proof}
Part a) follows from the compactness of $\operatorname{supp}\varphi$ and
$K$: for $(y,\sigma) \in K$, the function
$T_{1/\lambda,x_0} \circ T_{\sigma,y}(\varphi)$ is supported in
$B(x_0 + y/\lambda, \sigma/\lambda)$, which is contained in $\Omega$ for
$\lambda$ large enough.

For part b), note that for $\lambda > \bar{\lambda}_1$ (from part a)),
    \begin{equation}\label{energy decomposition}
        f_\lambda \circ T_{1/\lambda,x_0} \circ T_{\sigma,y}(\varphi) = \int_{\mathbb R^n}\left[|\nabla^m T_{\sigma,y}(\varphi)|^2 + \alpha(x)T_{\sigma,y}^2(\varphi)\right]dx + \int_\Omega \bar\alpha(x)\left[T_{1/\lambda,x_0} \circ T_{\sigma,y}(\varphi)\right]^2 dx.
    \end{equation}
The first term is $< S_{m,n} + \epsilon$ on $\partial K$ by
Corollary~\ref{test family degree}.  For the second term, since
$T_{1/\lambda,x_0} \circ T_{\sigma,y}(\varphi)$ is supported in
$B(x_0 + y/\lambda, \sigma/\lambda) \subset B(x_0, (r+\sigma_2)/\lambda)$,
    \[
        \int_\Omega \bar\alpha\left[T_{1/\lambda,x_0} \circ T_{\sigma,y}(\varphi)\right]^2 dx \leq \left(\int_{B(x_0,(r+\sigma_2)/\lambda)} \bar\alpha^{n/(2m)}\,dx\right)^{2m/n} \to 0 \quad (\lambda \to \infty).
    \]
Thus \eqref{support and energy}\,b) holds for $\lambda$ large enough.
\end{proof}

\subsection{Proof of Theorem~\ref{Main Theorem}}

\begin{proof}[Proof of Theorem~\ref{Main Theorem}]
Write
\[
    \kappa:=2^{2m/n}S_{m,n},
    \qquad
    \Psi_\lambda:=(\beta_\lambda,\gamma_\lambda),
\]
and, for the set $K$ supplied by Corollary~\ref{test family degree},
define the test family
\begin{equation}\label{rescaled test family}
    \Pi_\lambda(y,\sigma)
       :=T_{1/\lambda,x_0}\circ T_{\sigma,y}(\varphi).
\end{equation}
For $\lambda$ as in Lemma~\ref{rescaled test family bounds}, this is a
map from $K$ into $\mathcal N(\Omega)$.  The definitions of the rescalings give
\begin{equation}\label{test family coordinate identity}
    \Psi_\lambda\circ\Pi_\lambda=\Theta.
\end{equation}
In particular, \eqref{degree of test family} gives the nonzero degree
required in Theorem~\ref{abstract two level criterion}.

It is enough first to prove the assertion for
\[
    0<\epsilon<\min\{c(\alpha)-S_{m,n},
                         \kappa-S_{m,n}\}.
\]
Indeed, for an arbitrary prescribed $\epsilon>0$, one applies the
argument below with a smaller positive number satisfying this inequality
and smaller than $\epsilon$.  Choose $\varphi$ and $K$ as above, and
introduce the four levels
\begin{align*}
    a_\lambda&:=\inf\{f_\lambda(v):v\in \mathcal N(\Omega),
          \ \beta_\lambda(v)=0,\ \gamma_\lambda(v)\geq\tfrac13\},\\
    b_\lambda&:=\max_{\partial K}f_\lambda(\Pi_\lambda),\\
    d_\lambda&:=\inf\{f_\lambda(v):v\in \mathcal N(\Omega),
          \ \beta_\lambda(v)=0,\ \gamma_\lambda(v)=\tfrac13\},\\
    e_\lambda&:=\max_Kf_\lambda(\Pi_\lambda).
\end{align*}
Lemma~\ref{levels on the constraint} and
Lemma~\ref{rescaled test family bounds} yield
\begin{equation}\label{first four-level chain}
    S_{m,n}<a_\lambda,
    \qquad
    b_\lambda<S_{m,n}+\epsilon<c(\alpha)\leq d_\lambda.
\end{equation}
Moreover, $(0,\sigma_2)\in\partial K$ and radial symmetry give
\[
    \beta(T_{\sigma_2,0}\varphi)=0,
    \qquad
    \gamma(T_{\sigma_2,0}\varphi)>\tfrac13,
\]
so $a_\lambda\leq b_\lambda$.  Thus
\[
    S_{m,n}<a_\lambda\leq b_\lambda
       <\min\{d_\lambda,\kappa\}.
\]
By Corollary~\ref{constrained PS bridge}, $f_\lambda|_{\mathcal N(\Omega)}$
satisfies the Palais--Smale condition in $(S_{m,n},\kappa)$.
Theorem~\ref{abstract two level criterion}(i) therefore gives a
constrained critical point $v_{\lambda,\epsilon}$ satisfying
\begin{equation}\label{first critical level estimate}
    S_{m,n}<f_\lambda(v_{\lambda,\epsilon})
       <S_{m,n}+\epsilon.
\end{equation}

The Lagrange multiplier is $f_\lambda(v_{\lambda,\epsilon})$.  Hence
\begin{equation}\label{critical scaling}
    u_{\lambda,\epsilon}
       :=f_\lambda(v_{\lambda,\epsilon})^{1/(q-2)}
          v_{\lambda,\epsilon}
\end{equation}
is a nontrivial solution of the Euclidean equation, and therefore of
\eqref{Hyperbolic equation}; its energy is
\begin{equation}\label{scaled solution energy}
    I(u_{\lambda,\epsilon})
       =\frac mn f_\lambda(v_{\lambda,\epsilon})^{n/(2m)}.
\end{equation}
No positivity property of the perturbed operator has entered this
argument, and by Proposition~\ref{higher-order sign obstruction} none is
available under the present hypotheses.

The quantifiers in \eqref{first critical level estimate} are the
following: to each sufficiently small $\epsilon>0$ there corresponds a
$\Lambda(\epsilon)$ for which the estimate holds at every
$\lambda\geq\Lambda(\epsilon)$.  Choose $\epsilon_j\downarrow0$ and a
strictly increasing sequence $\Lambda_j\to\infty$ with
$\Lambda_j\geq\Lambda(\epsilon_j)$.  For
$\lambda\in[\Lambda_j,\Lambda_{j+1})$, select the solution constructed
with $\epsilon_j$.  The resulting family, denoted by $u_\lambda$, obeys
\begin{equation}\label{first family asymptotic}
    f_\lambda\left(
       \frac{u_\lambda}{\|u_\lambda\|_{L^q}}\right)
       \longrightarrow S_{m,n}
       \quad\text{as }\lambda\to\infty.
\end{equation}

We now assume
\[
    \|\alpha\|_{L^{n/(2m)}(\mathbb R^n)}
       <S_{m,n}(2^{2m/n}-1).
\]
The test function may then be chosen, in addition, so that
\begin{equation}\label{stronger test function}
    \int_{B(0,1)}|\nabla^m\varphi|^2\,dx
       <2^{2m/n}S_{m,n}
          -\|\alpha\|_{L^{n/(2m)}(\mathbb R^n)}.
\end{equation}
Indeed, the right-hand side is strictly larger than $S_{m,n}$.  Using
\eqref{energy decomposition}, H\"older's inequality, and the uniform
vanishing of the $\bar\alpha$ term gives, for all sufficiently large
$\lambda$,
\begin{equation}\label{energy below doubling}
    e_\lambda=\max_K f_\lambda(\Pi_\lambda)<\kappa.
\end{equation}
The nonzero degree in \eqref{degree of test family} also gives
$d_\lambda\leq e_\lambda$.  Combining this with
\eqref{first four-level chain}, we obtain the complete chain
\begin{equation}\label{complete four-level chain}
    S_{m,n}<a_\lambda\leq b_\lambda
       <d_\lambda\leq e_\lambda<\kappa.
\end{equation}
Part (ii) of Theorem~\ref{abstract two level criterion} gives a second
constrained critical point $\widehat v_{\lambda,\epsilon}$ such that
\[
    d_\lambda\leq
       f_\lambda(\widehat v_{\lambda,\epsilon})
       \leq e_\lambda.
\]
After the scaling \eqref{critical scaling}, it gives a second nontrivial
solution $\widehat u_{\lambda,\epsilon}$, and
\[
    f_\lambda(v_{\lambda,\epsilon})
       \leq b_\lambda<d_\lambda
       \leq f_\lambda(\widehat v_{\lambda,\epsilon}).
\]
Thus the two solutions are distinct, even modulo multiplication by $-1$.
To make the quantifiers explicit in the smallness case, choose
$\epsilon_j\downarrow0$ and test functions $\varphi_j$ satisfying both
\eqref{test function properties} and \eqref{stronger test function}.
Choose $\Lambda_j\uparrow\infty$ so that, for every
$\lambda\geq\Lambda_j$, all support and boundary estimates of
Lemma~\ref{rescaled test family bounds}, the estimate
$e_{\lambda,j}<\kappa$, and both min--max constructions hold for
$\varphi_j$ and its associated set $K_j$.  On each interval
$[\Lambda_j,\Lambda_{j+1})$, select the two critical points obtained
from $\varphi_j$.  This produces the two families asserted in the theorem
while retaining \eqref{first family asymptotic}.
\end{proof}

\medskip
\noindent\textbf{Acknowledgements.}
This work was supported by the National Natural Science Foundation of
China (No.\ 12571127).

\medskip
\noindent\textbf{Declaration of competing interest.}
The authors declare that they have no known competing financial interests
or personal relationships that could have appeared to influence the work
reported in this paper.

\medskip
\noindent\textbf{Data availability.}
No data were used for the research described in this article.

\end{document}